\documentclass[9pt]{article}
\usepackage{amsfonts}
\usepackage{amsmath}
\usepackage{amssymb}

\numberwithin{equation}{section}
\newtheorem{Theorem-}{Theorem}[section]
\newtheorem{Corollary-}{Corollary}[section]
\newtheorem{Definition-}{Definition}[section]
\newtheorem{Example-}{Example}[section]
\newtheorem{Lemma-}{Lemma}[section]
\newtheorem{Proposition-}{Proposition}[section]
\newtheorem{Remark-}{Remark}[section]
\newenvironment{Proof-}[1][Proof]{\noindent\textbf{#1.} }{\ \rule{0.5em}{0.5em}}

\newcommand{\ds}{\displaystyle}
\begin{document}
\title{\textbf{Generalized translation operator and the Heat equation for the canonical Fourier Bessel transform}}
\scrollmode
\author{Sami Ghazouani
\thanks{ Faculty of Sciences of Bizerte, (UR17ES21),"Dynamical Systems and their applications", University of Carthage, 7021 Jarzouna, Bizerte, Tunisia. E-mail addresses: Sami.Ghazouani@ipeib.rnu.tn}, Jihed Sahbani
\thanks{ Faculty of Sciences of Tunis, University of Tunis El Manar, Tunis, Tunisia}}
\date{ }

\maketitle \vspace{0mm}
\begin{abstract}
The aim of this paper is to introduce a translation operator associated to the canonical Fourier Bessel transform $\mathcal{F}_{\nu}^{\mathbf{m}}$ and study some of the important properties. We derive a convolution product for this transform and as application we study the heat equation and the heat semigroup related to $\Delta_{\nu}^{\mathbf{m}}.$
\end{abstract}
\noindent \textit{AMS CLASSIFICATION: 33D15; 47A05}
\section{Introduction}
\paragraph{}
The linear canonical transform (LCT) is a class of linear integral transform with matrix parameter $\mathbf{m}\in SL(2,\mathbb{R})$ \cite{Bultheel, Kutay,Wolf1}. It includes many well-known transforms such as the Fourier transform, the fractional Fourier transform, the Fresnel transform etc. These integral transforms are of importance in several areas of physics and mathematics \cite{Bultheel,Collins,Moshinsky,Wolf1}. As is well known, a wide number of papers have been successfully devoted to the extension of the theory of LCT to some other integral transforms such that the Hankel transform, the Dunkl transform \cite{De Jeu, Dunkl} and the Fourier Bessel transform \cite{Trimeche}. In \cite{Wolf}, Wolf enlarged the concept of LCT to the Hankel setting. He built a class of linear integral transform with a kernel involving a Bessel function and matrix parameter $\mathbf{m}\in SL(2,\mathbb{R}).$ Special cases of this
transformation are the Hankel transform and the fractional Hankel transform \cite{Kerr}. In \cite{Ghazouani}, the authors introduced the Dunkl linear canonical transform (DLCT) which is a generalization of the LCT in the framework of Dunkl transform \cite{De Jeu}. DLCT includes many well-known transforms such as the Dunkl transform \cite{De Jeu, Dunkl}, the fractional Dunkl transform \cite{Ghazouani1,Ghazouani2} and the canonical Fourier Bessel transform \cite{Dhaouadi, Ghazouani}. In \cite{Dhaouadi}, the authors developed a harmonic analysis related the canonical Fourier Bessel transform $\mathcal{F}_{\nu}^{\mathbf{m}}.$ Several properties, such as a Riemann-Lebesgue lemma, inversion formula, operational formulas, Plancherel
theorem and Babenko inequality and several uncertainty inequalities are established.
\paragraph{}
In the present work, we continue the analysis begun in \cite{Dhaouadi} by studying the translation operator and the convolution product related to this transformation. Following the framework of Delsarte \cite{Delsarte} and Levitan \cite{Levitan}, we introduce a generalized translation
$T_{x}^{\nu,\mathbf{m}}f(y)=u(x,y) \ (x,y\geq 0,\ \mathbf{m}\in SL(2,\mathbb{R}))$ of a function $f\in\mathcal{C}^{2}([0,+\infty[)$ as the solution to the following Cauchy problem
\begin{eqnarray*}
\left\{\begin{array}{l}
\Delta_{\nu,x}^{\mathbf{m}} \ u(x,y)=\Delta_{\nu,y}^{\mathbf{m}} \ u(x,y), \\
u(x,0)=f(x),  \ \frac{\partial}{\partial x}u(x,0)=0,\\
\end{array}
\right.
\end{eqnarray*}
where
$
\Delta_{\nu}^{\mathbf{m}}=\frac{d^{2}}{dx^{2}}+\left( \frac{2\nu +1}{x}-2i
\frac{d}{b}  x\right) \frac{d}{dx}-\left( \frac{d^{2}}{b^{2}}x^{2}+2i\left( \nu
+1\right) \frac{d}{b}\right)
$
and where $\Delta_{\nu,x}^{\mathbf{m}}$ act in the $x$ variable. We prove that the solution of the above problem can be written explicitly as
$$T_{x}^{\nu,\mathbf{m}}f(y)=u(x,y)=e^{\frac{i}{2}\frac{d}{b}(x^{2}+y^{2})}T_{x}^{\nu }\left[e^{-\frac{i}{2}\frac{d}{b}s^{2}}f(s)\right](y)
$$
where $T_{x}^{\nu}$ is the generalized translation operators associated with the Bessel operator $\Delta_{\nu}$ \cite{Levitan}. We also study some of the important properties of the translation operator $T_{x}^{\nu,\mathbf{m}}$ such that the continuity property with respect to the norm $\| .\|_{p,\nu}.$ Next, we give a generalized convolution product $\underset{\nu,\mathbf{m}}{*}$ on $(0,+\infty)$ tied to the differential operator $\Delta_{\nu}^{\mathbf{m}}$, by putting
\begin{eqnarray*}
f\underset{\nu,\mathbf{m}}{*}g(x)=\int_{0}^{+\infty }\left[T_{x}^{\nu,\mathbf{m}}f
(y)\right] \ \left[e^{-i\frac{d}{b}y^{2}}g(y)\right]y^{2\nu +1}dy
\end{eqnarray*}
and study some of its very basic properties. As application, we conclude this paper by studying the heat equation associated to $\sigma\Delta_{\nu}^{\mathbf{m}^{-1}}:$
\begin{eqnarray*}
\left\{\begin{array}{l}
\frac{\partial}{\partial t}u(t,x)=\sigma \Delta_{\nu}^{\mathbf{m}^{-1}}u(t,x), \ (t,x)\in (0,+\infty)\times\mathbb{R}  \\
u(0,x)=f(x).
\end{array}
\right.
\end{eqnarray*}
\paragraph{}
This paper is organized as follows. Section $2$ is devoted to an overview of the canonical Fourier Bessel transform. In section $3$ we introduce and study the translation operator and the convolution product related to the canonical Fourier Bessel transform $\mathcal{F}_{\nu}^{\mathbf{m}}$ and we establish some of its properties. In section $4$ we study the heat equation associated to $\Delta_{\nu}^{\mathbf{m}^{-1}}.$
\section{A brief survey of the canonical Fourier Bessel transform}
\paragraph{}
The aim of this section is to give a brief review of the theory of canonical Fourier Bessel transform that is relevant to the succeeding sections\cite{Dhaouadi}. Throughout this paper, $\nu$ denotes a real number such that $\nu>-\frac{1}{2}.$
\subsection{Notations and functional spaces}
We denote by:

\begin{quote}
$\bullet \ \mathcal{C}_{*,0}(\mathbb{R})$ the space of even continuous functions on $\mathbb{R}$ and vanishing at infinity. We provide $\mathcal{C}_{*,0}(\mathbb{R})$ with the topology of uniform convergence.\\
$\bullet \ \mathcal{C}_{*,c}(\mathbb{R})$ the space of even continuous functions on $\mathbb{R}$ and with compact support.\\
$\bullet \ \mathcal{C}_{*,b}(\mathbb{R})$ the space of even and bounded continuous functions on $\mathbb{R}.$ We provide $\mathcal{C}_{*,b}(\mathbb{R})$ with the topology of uniform convergence.\\
$\bullet \ \mathcal{L}_{p,\nu }$ the Lebesgue space of measurable functions on $[0,+\infty [ $ such that
\begin{eqnarray*}
\left\| f\right\| _{p,\nu }&=&\left(\ds\int_{0}^{+\infty }\left\vert f(y)\right\vert ^{p}y^{2\nu+1}dy\right)^{\frac{1}{p}}<\infty,  \  \mbox{if} \ 1\leq p<\infty, \\
\left\| f\right\| _{\infty}&=&\underset{y\in[0,\infty[}{\mbox{ess sup}}|f(y)|<\infty,  \  \mbox{if} \ p=\infty.
\end{eqnarray*}
We provide $\mathcal{L}_{p,\nu }$ with the topology defined by the norm $\left\Vert . \right\Vert _{p,\nu }$.\\
$ \bullet \ \mathcal{L}_{2,\nu }$ the Hilbert space equipped with the inner product $\langle .,.\rangle_{\nu}$ given by
$$\langle f,g\rangle_{\nu}=\int_{0}^{+\infty }f(y)\overline{g(y)}y^{2\nu+1}dy.
$$
$\bullet \ \mathcal{E}_{\ast}(\mathbb{R})$ the space of even $\mathcal{C}^{\infty }$-functions on $\mathbb{R}.$ We provide it with the topology of uniform
convergence on all compacts of $\mathbb{R},$ for functions and their derivatives.\\
$ \bullet \ \mathcal{S}_{\ast }(\mathbb{R})$ the Schwartz space of even $\mathcal{C}^{\infty }$-functions on $\mathbb{R}$ and rapidly decreasing together with their
derivatives. We provide $\mathcal{S}_{\ast }(\mathbb{R})$ with the topology defined by the seminorms:
\begin{equation*}
q_{n,m}(f)=\sup_{x\geq 0}(1+x^{2})^{n}\left\vert \frac{d^{m}}{dx^{m}}
f(x)\right\vert; \ f\in \mathcal{S}_{\ast }(\mathbb{R}), \ n,m\in \mathbb{N}.
\end{equation*}
$\bullet \ \mathcal{D}_{\ast}(\mathbb{R})$ the space of even $\mathcal{C}^{\infty}$-functions on $\mathbb{R}$ with compact support. We have
$$\mathcal{D}_{\ast}(\mathbb{R})=\underset{a\geq0}{\cup}\mathcal{D}_{\ast ,a}(\mathbb{R}),
$$
where $\mathcal{D}_{\ast ,a}(\mathbb{R})$ is the space of even $\mathcal{C}^{\infty}$-functions on $\mathbb{R}$ with support in the interval $[-a,a].$ We provide $\mathcal{D}_{\ast ,a}(\mathbb{R})$ with the topology of uniform convergence of functions and their derivatives.\\
$\bullet \ \mathbb{H}_{\ast}(\mathbb{C})$ the space of even entire functions on $\mathbb{C}$ rapidly decreasing and of exponential type. We have
$$\mathbb{H}_{\ast}(\mathbb{C})=\underset{a\geq0}{\cup}\mathbb{H}_{a}(\mathbb{C}),
$$
where $\mathbb{H}_{a}(\mathbb{C})$ is the space of even entire functions $f: \mathbb{C}\rightarrow \mathbb{C}$ and satisfying for all $m\in \mathbb{N}$
\begin{equation*}
P_{m}(f)=\sup_{\lambda \in \mathbb{C}}(1+|\lambda|)^{m}|f(\lambda)|e^{-a|\mathrm{Im}(\lambda )|}< \infty.
\end{equation*}
We provide $\mathbb{H}_{a}(\mathbb{C})$ with the topology defined by the seminorms $P_{m}, \ m\in \mathbb{N}.$
\end{quote}
\subsection{Canonical Fourier Bessel transform}
\paragraph{}
Throughout this paper, we denote by $\mathbf{m}=\left(\begin{array}{cc}a & b \\c & d\end{array}\right)$ an arbitrary matrix in $SL(2,\mathbb{R}).$ For $\mathbf{m}\in SL(2,\mathbb{R})$ such that $b\neq0,$ the canonical Fourier Bessel transform of a function $f\in \mathcal{L}_{1,\nu }$ is defined by \cite{Dhaouadi,Ghazouani}
\begin{eqnarray*}
\mathcal{F}_{\nu}^{\mathbf{m}}f(x)=\ds\frac{c_{\nu}}{(ib)^{\nu+1}} \ds\int_{0}^{+\infty}\ K_{\nu}^{\mathbf{m}}(x,y) f(y)y^{2\nu+1}dy,
\end{eqnarray*}
where $c_{\nu}^{-1}=2^{\nu}\Gamma(\nu+1)$ and
\begin{eqnarray*}
K_{\nu}^{\mathbf{m}}(x,y)=e^{\frac{i}{2}(\frac{d}{b}x^2+\frac{a}{b}y^{2})
} \ j_{\nu}\left(\frac{xy}{b}\right).
\end{eqnarray*}
Here $j_{\nu }$ denotes the normalized Bessel function of order $\nu>-1/2$ and defined by \cite{George,Watson}
\begin{equation*}
j_{\nu }(x)=2^{\nu}\Gamma(\nu+1)\ds\frac{J_{\nu}(x)}{x^{\nu}}=\sum_{n=0}^{+\infty }(-1)^{n}\frac{\Gamma (\nu +1)}{n!\Gamma (n+\nu
+1)}\left(\frac{x}{2}\right)^{2n}; \ x\in \mathbb{C}.
\end{equation*}
It is well known that, for every $y\in \mathbb{C},$ the function $j_{\nu }(y\cdot):x\mapsto j_{\nu }(y x)$
is the unique solution of
\begin{eqnarray*}
\left\{\begin{array}{ll}
\Delta _{\nu } f=-y ^{2}f \\
f(0)=1,  \ f^{'}(0)=0,
\end{array}\right.
\end{eqnarray*}
where $\Delta _{\nu }$ is the Bessel operator given by
\begin{equation*}
\Delta _{\nu }=\frac{d^{2}}{dx^{2}}+\frac{2\nu +1}{x}\frac{d}{dx}.
\end{equation*}
For every $\nu> -\frac{1}{2},$ the normalized Bessel function $j_{\nu}$ has the Mehler integral representation
\begin{equation*}
j_{\nu}(x)=
\ds\frac{2\Gamma(\nu+1)}{\sqrt{\pi}\Gamma(\nu+1/2)}\ds\int_{0}^{1}\left(1-t^{2}\right)^{\nu-\frac{1}{2}} \cos(xt) dt.
\end{equation*}
Consequently, we have
\begin{eqnarray*}
\forall x \in \mathbb{C}, \quad |j_{\nu}(x)|\leq e^{|\mathrm{Im}(x)|}.
\end{eqnarray*}
The normalized Bessel function $j_{\nu}$ satisfies, for all $x,z \in [0,+\infty[$, the following product formula \cite{Watson}
\begin{equation*}
j_{\nu}(xy)j_{\nu}(yz)=
\ds\frac{\Gamma(\nu+1)}{\sqrt{\pi}\Gamma(\nu+1/2)}\ds\int_{0}^{\pi}j_{\nu}\left(y\sqrt{x^{2}+z^{2}-2xz\cos(\theta)}\right) \sin^{2\nu}(\theta)d\theta.
\end{equation*}
\paragraph{}
It is easy to see that, for each $y\in \mathbb{R},$ the kernel $K_{\nu}^{\mathbf{m}}(\cdot,y)$ of the canonical Fourier Bessel transform $\mathcal{F}_{\nu}^{\mathbf{m}}$ is the unique solution of \cite{Dhaouadi}:
\begin{eqnarray*}
\left\{\begin{array}{ll}
\Delta_{\nu}^{\mathbf{m}} \ K_{\nu}^{\mathbf{m}}(\cdot,y)=-\frac{y^{2}}{b^{2}} \ K_{\nu}^{\mathbf{m}}(\cdot,y) , &  \\
K_{\nu}^{\mathbf{m}}(0,y)=e^{\frac{i}{2}\frac{a}{b}y^{2}},  \ \frac{d}{dx}K_{\nu}^{\mathbf{m}}(0,y)=0,
\end{array}\right.
\end{eqnarray*}
where for each $\mathbf{m}\in SL(2,\mathbb{R})$ such that $b\neq0,\ \Delta_{\nu}^{\mathbf{m}}$ denotes the differential operator
\begin{eqnarray}
\Delta_{\nu}^{\mathbf{m}}=\frac{d^{2}}{dx^{2}}+\left( \frac{2\nu +1}{x}-2i
\frac{d}{b}  x\right) \frac{d}{dx}-\left( \frac{d^{2}}{b^{2}}x^{2}+2i\left( \nu
+1\right) \frac{d}{b}\right).\label{equation2}
\end{eqnarray}
\paragraph{}
Many properties of the Fourier Bessel transform carry over to the canonical Fourier Bessel transform. In particular:
\begin{Theorem-} \cite{Dhaouadi} \mbox{}\\
$(1)$ Let $\mathbf{m}\in SL(2,\mathbb{R}).$\\
$(a)$\textbf{The reversibility property}: For all $f \in \mathcal{L}_{1,\nu }$ with $\mathcal{F}_{\nu}^{\mathbf{m}}
f \in \mathcal{L}_{1,\nu }$,
\begin{eqnarray*}
\left(\mathcal{F}_{\nu}^{\mathbf{m}}\circ \mathcal{F}_{\nu}^{\mathbf{m}^{-1}}\right)f=\left(\mathcal{F}_{\nu}^{\mathbf{m}^{-1}}\circ\mathcal{F}_{\nu}^{\mathbf{m}}\right) f=f, \  a. e.
\end{eqnarray*}
$(b)$ The canonical Fourier Bessel transform $\mathcal{F}_{\nu}^{\mathbf{m}}$ is a topological isomorphism from $\mathcal{S}_{\ast }(\mathbb{R})$ into itself. More precisely, for each $(n,m)\in \mathbb{N}\times\mathbb{N}$ there exist $n_{1},\dots, n_{k}, m_{1},\dots, m_{k} \in \mathbb{N}$ and $c>0$ such that
\begin{eqnarray}
q_{n,m}\left(\mathcal{F}_{\nu}^{\mathbf{m}}f\right)\leq c\ds\sum_{i=1}^{k}q_{n_{i},m_{i}}(f).\label{equat3}
\end{eqnarray}
$(2)$\textbf{Operational formulas}: Let $\mathbf{m}\in SL(2,\mathbb{R})$ such that $b\neq 0$ and $f\in\mathcal{S}_{\ast }(\mathbb{R}).$ Then\\
$(a)$
\begin{eqnarray*}
\mathcal{F}_{\nu}^{\mathbf{m}}\left[y^{2}f(y)\right](x)=-b^{2}\Delta_{\nu}^{\mathbf{m}}\left[\mathcal{F}_{\nu}^{\mathbf{m}}f\right](x).
\end{eqnarray*}
$(b)$
\begin{eqnarray}
x^{2}\mathcal{F}_{\nu}^{\mathbf{m}}(f)(x)=-b^{2}\mathcal{F}_{\nu}^{\mathbf{m}}\left[\Delta_{\nu}^{\mathbf{m}^{-1}}f\right](x).\label{equa1}
\end{eqnarray}
$(3)$\textbf{Babenko inequality}: Let $\mathbf{m}\in SL(2,\mathbb{R})$ such that $b\neq0.$ Let $p$ and $q$ be real numbers such that $1<p\leq 2$ and $\frac{1}{p}+\frac{1}{q}=1.$ Then $\mathcal{F}_{\nu}^{\mathbf{m}}$ extends to a bounded linear operator
on $\mathcal{L}_{p,\nu }$ and we have
\begin{equation}
\left\Vert \mathcal{F}_{\nu}^{\mathbf{m}}f\right\Vert _{q,\nu }\leq  |b|^{(\nu+1)(\frac{2}{q}-1)}\left( \frac{(c_{\nu}p)^{\frac{
1}{p}}}{(c_{\nu}q)^{\frac{1}{q}}}\right) ^{\nu +1} \| f\|_{p,\nu}.\label{babenko}
\end{equation}
\label{theorem1}
\end{Theorem-}
\section{Translation operators associated with $\Delta_{\nu}^{\mathbf{m}}$}
\paragraph{}
We shall begin the section by recalling some well-known facts about the generalized translation operators $T_{x}^{\nu }$ associated
with the Bessel operator $\Delta_{\nu}.$ Next we introduce and study a translation operators associated with $\Delta_{\nu}^{\mathbf{m}}.$
\subsection{Generalized translation operators associated with $\Delta_{\nu}$}
\paragraph{}
Levitan \cite{Levitan} defined the generalized translation operators $T_{x}^{\nu }$ associated
with $\Delta_{\nu}$ by:
\begin{eqnarray*}
T_{x}^{\nu }f(y)=u(x,y)
\end{eqnarray*}
where $u(x,y)$ is the unique solution of the Cauchy problem \cite{Levitan}:
\begin{eqnarray}
\left\{\begin{array}{l}
\left(\frac{\partial^{2}}{\partial x^{2}}+\frac{2\nu+1}{x}\frac{\partial}{\partial x}\right)u(x,y)=\left(\frac{\partial^{2}}{\partial y^{2}}+\frac{2\nu+1}{y}\frac{\partial}{\partial y}\right)u(x,y), \\
u(x,0)=f(x), \ f\in \mathcal{E}_{*}(\mathbb{R}), \\
\frac{\partial}{\partial x}u(x,0)=0.
\end{array}
\right.\label{prob}
\end{eqnarray}
The solution of (\ref{prob}) can be represented in the following form \cite{Levitan}:
\begin{eqnarray*}
T_{x}^{\nu }f(y)=\frac{\Gamma(\nu+1)}{\sqrt{\pi}\Gamma(\nu+1/2)}\ds\int_{0}^{\pi}f\left(\sqrt{x^{2}+y^{2}-2xy\cos(\theta)}\right)\left(\sin(\theta)\right)^{2\nu} d\theta.
\end{eqnarray*}
It is interesting to note that $T_{x}^{\nu }f(y)$ is an even function with respect to both variables and then we can assume that $x\geq 0$ and $y\geq0.$ By a
change of variables we get
\begin{equation}
T_{x}^{\nu }f(y)=\int_{0}^{+\infty }f(z) \ \mathcal{W}_{\nu}(x,y,z) \ z^{2\nu +1}dz; \ x,y\geq0.\label{equat1}
\end{equation}
Here, the kernel $\mathcal{W}_{\nu}(x,y,z)$ is given by
\begin{equation*}
\mathcal{W}_{\nu}(x,y,z)=\frac{2^{1-2\nu}\Gamma(\nu+1)}{\sqrt{\pi}\Gamma(\nu+1/2)}\frac{\Delta^{2\nu-1}}{(xyz)^{2\nu}}\upharpoonleft_{[|x-y|, \ x+y]}
\end{equation*}
where
\begin{equation*}
\Delta=\frac{1}{4}\sqrt{(x+y+z)(x+y-z)(x-y+z)(y+z-x)}
\end{equation*}
denotes the area of the triangle with sides $x,y,z>0$ and $\upharpoonleft_{A}$ is the indicator function of $A.$
We list some standard properties of the generalized translation operators $T_{x}^{\nu}$ that we shall use in this article.
\begin{Theorem-} \cite{Levitan} The operators $T_{x}^{\nu} ,x\in \mathbb{R},$ satisfy:\\
$1)$ \textbf{Symmetry:} $ T_{0}^{\nu}=\mbox{identity}$ and $T_{x}^{\nu}f(y)=T_{y}^{\nu}f(x).$\\
$2)$ \textbf{Linearity:} $T_{x}^{\nu }\left[af+g\right](y)=aT_{x}^{\nu }f(y)+T_{x}^{\nu }g(y).$\\
$3)$ \textbf{Positivity:} If $f(y)\geq 0$ then $T_{x}^{\nu}f(y)\geq0.$\\
$4)$ \textbf{Compact support:} If $f(y)=0$ for $y\geq a$ then $T_{x}^{\nu}f(y)=0$ for $|x-y|\geq a.$\\
$5)$ \textbf{Product formula:} For all $x,y \in \mathbb{R}$ and $\lambda \in \mathbb{C},$ we have the product formula
$$T_{x}^{\nu }\left[s\mapsto j_{\nu}(\lambda s)\right](y)=j_{\nu}(\lambda x) j_{\nu}(\lambda y).
$$
$5)$ \textbf{Continuity of translation operators $T_{x}^{\nu}$:} The operator $T_{x}^{\nu}$ is continuous from:\\
$\bullet$ $\mathcal{C}(\mathbb{R})$ into itself. (The space $\mathcal{C}(\mathbb{R})$ of all continuous functions on $\mathbb{R}$ is equipped with the topology of uniform convergence
on compact sets)\\
$\bullet$ $\mathcal{E}_{*}(\mathbb{R})$ into itself. \\
$\bullet$ $\mathcal{L}_{p,\nu }$ into itself. More precisely, for all $f\in \mathcal{L}_{p,\nu }, \ p\in [1,+\infty],$ and $x\geq 0,$ the function $T_{x}^{\nu }f$ is defined
almost everywhere on $[0,+\infty[,$ belongs to $\mathcal{L}_{p,\nu }$ and we have
$$\left\|T_{x}^{\nu}f\right\|_{p,\nu}\leq \left\|f\right\|_{p,\nu}.
$$
$6)$ \textbf{Commutativity:}\\
$\bullet$ The following equality hold in $\mathcal{C}(\mathbb{R}): \ T_{x}^{\nu}\circ T_{y}^{\nu}=T_{y}^{\nu}\circ T_{x}^{\nu}.$\\
$\bullet$ The following equality hold in $\mathcal{E}_{*}(\mathbb{R}): \ \Delta_{\nu}\circ T_{x}^{\nu}=T_{x}^{\nu}\circ\Delta_{\nu}.$\\
$7)$ \textbf{Self-adjointness of $T_{x}^{\nu}$:} Let $f\in \mathcal{L}_{1,\nu }$ and $g\in \mathcal{C}_{*,b}(\mathbb{R}).$ Then
$$\ds\int_{0}^{+\infty}T_{x}^{\nu}f(y)g(y) y^{2\nu+1} dy=\ds\int_{0}^{+\infty}f(y)T_{x}^{\nu}g(y)  y^{2\nu+1} dy.
$$
In particular:\\
$\bullet$ For $g(y)=1,$
$$\ds\int_{0}^{+\infty}T_{x}^{\nu}f(y)y^{2\nu+1} dy=\ds\int_{0}^{+\infty}f(y)y^{2\nu+1} dy.
$$
$\bullet$ For $g(y)=j_{\nu}(\lambda y), \ \lambda\in \mathbb{R},$
\begin{eqnarray*}
\mathcal{F}_{\nu}\left[T_{x}^{\nu}f\right](\lambda)=j_{\nu}(\lambda x)
\ \mathcal{F}_{\nu}f(\lambda).
\end{eqnarray*}
\end{Theorem-}
\subsection{Generalized translation operators associated with $\Delta_{\nu}^{\mathbf{m}}$}
We begin with the following Proposition
\begin{Proposition-} Let $\mathbf{m}\in SL(2,\mathbb{R})$ such that $b\neq0.$ For each $f\in \mathcal{E}_{*}(\mathbb{R})$ the problem
\begin{eqnarray}
\left\{\begin{array}{l}
\Delta_{\nu,x}^{\mathbf{m}} \ u(x,y)=\Delta_{\nu,y}^{\mathbf{m}} \ u(x,y), \\
u(x,0)=f(x),  \\
\frac{\partial}{\partial x}u(x,0)=0,
\end{array}
\right.\label{prob1}
\end{eqnarray}
has a unique solution given by:
\begin{eqnarray*}
u(x,y)=e^{\frac{i}{2}\frac{d}{b}(x^{2}+y^{2})}T_{x}^{\nu }\left[e^{-\frac{i}{2}\frac{d}{b}s^{2}}f(s)\right](y).
\end{eqnarray*}
Here $\Delta_{\nu,x}^{\mathbf{m}}$ act in the $x$ variable and given by (\ref{equation2}).
\label{translation}
\end{Proposition-}
\begin{Proof-} From the transmutation property
\begin{eqnarray*}
e^{-\frac{i}{2}\frac{d}{b}x^{2}}\circ \Delta_{\nu}^{\mathbf{m}}\ \circ e^{\frac{i}{2}\frac{d}{b}x^{2}}=\Delta_{\nu},
\end{eqnarray*}
the system (\ref{prob1}) is equivalent to
\begin{eqnarray*}
\left\{\begin{array}{l}
\left(\frac{\partial^{2}}{\partial x^{2}}+\frac{2\nu+1}{x}\frac{\partial}{\partial x}\right)\tilde{u}(x,y)=\left(\frac{\partial^{2}}{\partial y^{2}}+\frac{2\nu+1}{y}\frac{\partial}{\partial y}\right)\tilde{u}(x,y), \\
\tilde{u}(x,0)=e^{-\frac{i}{2}\frac{d}{b}x^{2}}f(x), \\
\frac{\partial}{\partial x}\tilde{u}(x,0)=0,
\end{array}
\right.
\end{eqnarray*}
where $\tilde{u}(x,y)=e^{-\frac{i}{2}\frac{d}{b}(x^{2}+y^{2})}u(x,y).$ By (\ref{prob}) we obtain
$$\tilde{u}(x,y)=T_{x}^{\nu }\left[e^{-\frac{i}{2}\frac{d}{b}s^{2}}f(s)\right](y)
$$
and then
\begin{eqnarray*}
u(x,y)=e^{\frac{i}{2}\frac{d}{b}(x^{2}+y^{2})}T_{x}^{\nu }\left[e^{-\frac{i}{2}\frac{d}{b}s^{2}}f(s)\right](y).
\end{eqnarray*}
\end{Proof-}
\begin{Definition-} Let $\mathbf{m}\in SL(2,\mathbb{R})$ such that $b\neq0.$ For $f\in \mathcal{C}(\mathbb{R}),$ we define the generalized translation operators associated with the operator $\Delta_{\nu}^{\mathbf{m}}$ by:
\begin{eqnarray*}
T_{x}^{\nu,\mathbf{m}}f(y)=e^{\frac{i}{2}\frac{d}{b}(x^{2}+y^{2})}T_{x}^{\nu }\left[e^{-\frac{i}{2}\frac{d}{b}s^{2}}f(s)\right](y),
\end{eqnarray*}
where $T_{x}^{\nu}$ is the generalized translation operator associated with $\Delta_{\nu}.$
\end{Definition-}
\paragraph{}
The following definition will be important in the sequel:
\begin{Definition-} \cite{Dhaouadi} The chirp multiplication operator $\mathbf{L}_{a }$ and the dilatation operator $\mathbf{D}_{a}$ are defined, respectively, by
\begin{eqnarray*}
\mathbf{L}_{a }f(x)&=&e^{\frac{ia }{2}x^{2}}f(x); \ a \in \mathbb{R}.\\ \mathbf{D}_{a }f(x)&=&\frac{1}{|a|^{\nu+1}}f\left(\frac{x}{a}\right); \ a \in \mathbb{R}\backslash\{0\}.
\end{eqnarray*}
In the case where $f$ is an even function, we have the following result:
\begin{eqnarray*}
\mathbf{D}_{a }f(x)&=&\mathbf{D}_{|a|}f(x); \ a \in \mathbb{R}\backslash\{0\}.
\end{eqnarray*}
The inverse of $\mathbf{L}_{a }$ and the inverse of $\mathbf{D}_{a}$ are given, respectively, by
\begin{eqnarray*}
\left(\mathbf{L}_{a }\right)^{-1}=\mathbf{L}_{-a}, \quad \left(\mathbf{D}_{a }\right)^{-1}=\mathbf{D}_{\frac{1}{a}}.
\end{eqnarray*}
\end{Definition-}
\paragraph{}
We give in the following proposition some useful properties of $\mathbf{L}_{a}$ and $\mathbf{D}_{a}.$
\begin{Proposition-}\cite{Dhaouadi} Let $1 \leq p <\infty$ and $a \in \mathbb{R}.$\\
$1)i)$ The operator $\mathbf{L}_{a }$ is an isometric isomorphism from:\\
$\bullet$ $ \mathcal{L}_{p,\nu}$ into itself.\\
$\bullet$ $ \mathcal{C}_{*,0}(\mathbb{R})$ into itself.\\
$ii)$ $\mathbf{L}_{a }$ is a unitary operator from $\mathcal{L}_{2,\nu}$ into itself.\\
$iii)$ The operator $\mathbf{L}_{a}$ is a topological isomorphism from:\\
$\bullet$ $\mathcal{S}_{\ast}(\mathbb{R})$ into itself.\\
$\bullet$ $\mathcal{D}_{\ast,r}(\mathbb{R})$ into itself.\\
$2)i)$ For each $a>0,$ the operator $\mathbf{D}_{a}$ is a topological isomorphism from $\mathcal{L}_{p,\nu}$ into itself and we have
$$\|\mathbf{D}_{a}f\|_{p,\nu}=a^{(\nu+1)(\frac{2}{p}-1)}\|f\|_{p,\nu}\ ; \ f\in \mathcal{L}_{p,\nu}.
$$
In particular, $\mathbf{D}_{a}$ is a unitary operator from $\mathcal{L}_{2,\nu}$ into itself. \\
$ii)$ For each $a\neq 0,$ the operator $\mathbf{D}_{a}$ is a topological isomorphism from $\mathcal{C}_{*,0}(\mathbb{R})$ into itself and we have
$$\|\mathbf{D}_{a}f\|_{\infty}=|a|^{-(\nu+1)}\|f\|_{\infty}\ ; \ f\in \mathcal{C}_{*,0}(\mathbb{R}).
$$
$iii)$ The operator $\mathbf{D}_{a}$ is a topological isomorphism from:\\
$\bullet$ $\mathcal{S}_{\ast}(\mathbb{R})$ into itself.\\
$\bullet$ $\mathcal{D}_{\ast,r}(\mathbb{R})$ into $\mathcal{D}_{\ast,ar}(\mathbb{R})$.
\end{Proposition-}
\paragraph{}
Now we are in a position to give some properties of the $T_{x}^{\nu,\mathbf{m}}.$
\begin{Proposition-} Let $\mathbf{m}\in SL(2,\mathbb{R})$ such that $b\neq0.$ The operators $T_{x}^{\nu,\mathbf{m}} ,x\in \mathbb{R},$ satisfy:\\
$1)$ \textbf{Symmetry:} $T_{0}^{\nu,\mathbf{m}}=\mbox{identity}$ and $T_{x}^{\nu,\mathbf{m}}f(y)=T_{y}^{\nu,\mathbf{m}}f(x).$\\
$2)$ \textbf{Linearity:} $T_{x}^{\nu,\mathbf{m}}\left[\lambda f+g\right](y)=\lambda T_{x}^{\nu,\mathbf{m}}f(y)+T_{x}^{\nu,\mathbf{m}}g(y).$\\
$3)$ \textbf{Compact support:} If $f(y)=0$ for $y\geq r$ then $T_{x}^{\nu,\mathbf{m}}f(y)=0$ for $|x-y|\geq r.$\\
$4)$ \textbf{Product formula:} For all $x, \ y$ and $ z \in \mathbb{R},$ we have the product formula
\begin{eqnarray*}
\begin{array}{l}
T_{x}^{\nu,\mathbf{m}}\left[K_{\nu}^{\mathbf{m}}(.,y)\right](z)=e^{-\frac{i}{2}\frac{a}{b}y^{2}}K_{\nu}^{\mathbf{m}}(x,y) K_{\nu}^{\mathbf{m}}(z,y).
\end{array}
\end{eqnarray*}
$5)$ For all $x,y>0$ we have
\begin{eqnarray*}
T_{x}^{\nu,\mathbf{m}}f(y)=\int_{0}^{+\infty }e^{-i\frac{d}{b}z^{2}}f(z) \ \mathcal{W}_{\nu}^{\mathbf{m}}(x,y,z) \ z^{2\nu +1}dz,
\end{eqnarray*}
where
\begin{eqnarray*}
\mathcal{W}_{\nu}^{\mathbf{m}}(x,y,z)=e^{\frac{i}{2}\frac{d}{b}(x^{2}+y^{2}+z^{2})}\mathcal{W}_{\nu}(x,y,z).
\end{eqnarray*}
$6)$ \textbf{Continuity of translation operators $T_{x}^{\nu,\mathbf{m}}$:} The operator $T_{x}^{\nu,\mathbf{m}}$ is continuous from:\\
$\bullet$ $\mathcal{C}(\mathbb{R})$ into itself. (The space $\mathcal{C}(\mathbb{R})$ of all continuous functions on $\mathbb{R}$ is equipped with the topology of uniform convergence
on compact sets)\\
$\bullet$ $\mathcal{E}_{*}(\mathbb{R})$ into itself. \\
$\bullet$ $\mathcal{L}_{p,\nu }$ into itself. More precisely, for all $f\in \mathcal{L}_{p,\nu }, \ p\in [1,+\infty],$ and $x\geq 0,$ the function $T_{x}^{\nu,\mathbf{m}}f$ is defined
almost everywhere on $[0,+\infty[,$ belongs to $\mathcal{L}_{p,\nu }$ and we have
\begin{eqnarray}
\left\|T_{x}^{\nu,\mathbf{m}}f\right\|_{p,\nu}\leq \left\|f\right\|_{p,\nu}.\label{inequa}
\end{eqnarray}
$7)$ \textbf{Commutativity:}\\
$\bullet$ The following equality hold in $\mathcal{C}(\mathbb{R}): \ T_{x}^{\nu,\mathbf{m}}\circ T_{y}^{\nu,\mathbf{m}}=T_{y}^{\nu,\mathbf{m}}\circ T_{x}^{\nu,\mathbf{m}}.$\\
$\bullet$ The following equality hold in $\mathcal{E}_{*}(\mathbb{R}): \ \Delta_{\nu}^{\mathbf{m}}\circ T_{x}^{\nu,\mathbf{m}}=T_{x}^{\nu,\mathbf{m}}\circ\Delta_{\nu}^{\mathbf{m}}.$\\
$8)$ Let $f\in \mathcal{L}_{1,\nu }$ and $g\in \mathcal{C}_{b}(\mathbb{R}).$ Then
\begin{eqnarray}
\int_{0}^{+\infty }\left[T_{x}^{\nu,\mathbf{m}}f
(y)\right] \ \left[e^{-i\frac{d}{b}y^{2}}g(y)\right]y^{2\nu +1}dy=\ds\int_{0}^{+\infty}\left[e^{-i\frac{d}{b}y^{2}}f(y)\right] \left[T_{x}^{\nu,\mathbf{m}}g
(y)\right]  y^{2\nu+1} dy.\label{equat}
\end{eqnarray}
$9)$ For all $f\in \mathcal{L}_{1,\nu },$ we have
\begin{eqnarray}
\mathcal{F}_{\nu}^{\mathbf{m}}\left[T_{x}^{\nu,\mathbf{m}^{-1}}f\right](\lambda)=e^{\frac{i}{2}\frac{d}{b}\lambda^{2}}\ K_{\nu}^{\mathbf{m}^{-1}}(x,\lambda)\mathcal{F}_{\nu}^{\mathbf{m}}f(\lambda).\label{fourier}
\end{eqnarray}
$10)$ For all $f\in \mathcal{L}_{p,\nu} , \ p\in ]1,2],$ we have
\begin{eqnarray*}
\mathcal{F}_{\nu}^{\mathbf{m}}\left[T_{x}^{\nu,\mathbf{m}^{-1}}f\right](\lambda)=e^{\frac{i}{2}\frac{d}{b}\lambda^{2}}\ K_{\nu}^{\mathbf{m}^{-1}}(x,\lambda)\mathcal{F}_{\nu}^{\mathbf{m}}f(\lambda), \ \mbox{a.e.}
\end{eqnarray*}
\end{Proposition-}
\begin{Proof-} Statements $1), \ 2)$ and $ 3)$ are obvious.\\
$4)$ From the definition of $T_{x}^{\nu,\mathbf{m}}$ it is easy to see that
\begin{eqnarray*}
T_{x}^{\nu,\mathbf{m}}\left[K_{\nu}^{\mathbf{m}}(\cdot,y)\right](z)&=&e^{\frac{i}{2}\frac{d}{b}(x^{2}+z^{2})}
T_{x}^{\nu}\left[s\mapsto e^{\frac{i}{2}\frac{a}{b}y^{2}}j_{\nu}(sy/b)\right](z)\\
&=&e^{\frac{i}{2}\frac{d}{b}(x^{2}+z^{2})}e^{\frac{i}{2}\frac{a}{b}y^{2}}T_{x}^{\nu}\left[s\mapsto j_{\nu}(sy/b)\right](z)\\
&=&e^{\frac{i}{2}\frac{d}{b}(x^{2}+z^{2})}e^{\frac{i}{2}\frac{a}{b}y^{2}}j_{\nu}(xy/b)j_{\nu}(zy/b)\\
&=&e^{-\frac{i}{2}\frac{a}{b}y^{2}}K_{\nu}^{\mathbf{m}}(x,y) K_{\nu}^{\mathbf{m}}(z,y).
\end{eqnarray*}
$5)$ This is a direct consequence of (\ref{equat1}).\\
$6)$ This follows from
$$T_{x}^{\nu,\mathbf{m}}f(y)=\left[\mathbf{L}_{\frac{d}{b},x}\circ \mathbf{L}_{\frac{d}{b},y}\circ T_{x}^{\nu}\circ \mathbf{L}_{-\frac{d}{b}}\right]f(y)
$$
and the fact that $\mathbf{L}_{\frac{d}{b}}, \   \mathbf{L}_{-\frac{d}{b}}, T_{x}^{\nu}$ are continuous from $\mathcal{C}(\mathbb{R})$ into itself, $\mathcal{E}_{*}(\mathbb{R})$ into itself and $\mathcal{L}_{p,\nu }$ into itself respectively.
 Now let $f\in\mathcal{L}_{p,\nu}.$ The function $T_{x}^{\nu,\mathbf{m}}f$ belongs to $\mathcal{L}_{p,\nu}$ and we have
\begin{eqnarray*}
\left\| T_{x}^{\nu,\mathbf{m}}f\right\|_{p,\nu}
&= &  \left\| T_{x}^{\nu}\left[\mathbf{L}_{-\frac{d}{b}}f\right]\right\|_{p,\nu}\\
&\leq&  \left\| \mathbf{L}_{-\frac{d}{b}}f\right\|_{p,\nu}=\left\| f\right\|_{p,\nu}.
\end{eqnarray*}
$7)$ $\bullet$ Let $f\in \mathcal{C}(\mathbb{R}).$ Then
 \begin{eqnarray*}
\left[T_{x}^{\nu,\mathbf{m}}\circ T_{y}^{\nu,\mathbf{m}}\right]f(z)&=&e^{\frac{i}{2}\frac{d}{b}(x^{2}+y^{2}+z^{2})}\left[T_{x}^{\nu}\circ T_{y}^{\nu}\right]\left[e^{-\frac{i}{2}\frac{d}{b}s^{2}}f(s)\right](z)\\
&=&e^{\frac{i}{2}\frac{d}{b}(x^{2}+y^{2}+z^{2})}\left[T_{y}^{\nu}\circ T_{x}^{\nu}\right]\left[e^{-\frac{i}{2}\frac{d}{b}s^{2}}f(s)\right](z)\\
&=&\left[T_{y}^{\nu,\mathbf{m}}\circ T_{x}^{\nu,\mathbf{m}}\right]f(z).
\end{eqnarray*}
$\bullet$ Let $f\in \mathcal{E}_{*}(\mathbb{R}).$ Then
 \begin{eqnarray*}
\left[\Delta_{\nu}^{\mathbf{m}}\circ T_{x}^{\nu,\mathbf{m}}\right]f(y)&=&e^{\frac{i}{2}\frac{d}{b}(x^{2}+y^{2})}\left[\Delta_{\nu}\circ T_{x}^{\nu}\right]\left[e^{-\frac{i}{2}\frac{d}{b}s^{2}}f(s)\right](y)\\
&=&e^{\frac{i}{2}\frac{d}{b}(x^{2}+y^{2})}\left[T_{x}^{\nu} \circ  \Delta_{\nu}\right]\left[e^{-\frac{i}{2}\frac{d}{b}s^{2}}f(s)\right](y)\\
&=&\left[T_{x}^{\nu,\mathbf{m}}\circ\Delta_{\nu}^{\mathbf{m}}\right]f(y).
\end{eqnarray*}
$8)$Let $f\in \mathcal{L}_{1,\nu }$ and $g\in \mathcal{C}_{*,b}(\mathbb{R}).$ Then
\begin{eqnarray*}
\int_{0}^{+\infty }\left[T_{x}^{\nu,\mathbf{m}}f
(y)\right] \ \left[e^{-i\frac{d}{b}y^{2}}g(y)\right]y^{2\nu +1}dy&=&
e^{\frac{i}{2}\frac{d}{b}x^{2}}\ds\int_{0}^{+\infty}T_{x}^{\nu }\left[e^{-\frac{i}{2}\frac{d}{b}s^{2}}f(s)\right](y)\left[e^{-\frac{i}{2}\frac{d}{b}y^{2}}g(y)\right] y^{2\nu+1} dy\\
&=&e^{\frac{i}{2}\frac{d}{b}x^{2}}\ds\int_{0}^{+\infty}\left[e^{-\frac{i}{2}\frac{d}{b}y^{2}}f(y)\right]T_{x}^{\nu }\left[e^{-\frac{i}{2}\frac{d}{b}s^{2}}g(s)\right](y) y^{2\nu+1} dy\\
&=&\ds\int_{0}^{+\infty}\left[e^{-i\frac{d}{b}y^{2}}f(y)\right] \left[T_{x}^{\nu,\mathbf{m}}g
(y)\right]  y^{2\nu+1} dy.
\end{eqnarray*}
$9)$Let $f\in \mathcal{L}_{1,\nu }.$ Then
\begin{eqnarray*}
\left[(ib)^{\nu+1}/c_{\nu}\right]\mathcal{F}_{\nu}^{\mathbf{m}}\left[T_{x}^{\nu,\mathbf{m}^{-1}}f\right](\lambda)&=&e^{\frac{i}{2}(\frac{d}{b}\lambda^{2}-\frac{a}{b}x^{2})}
\ds\int_{0}^{+\infty}T_{x}^{\nu }\left[e^{\frac{i}{2}\frac{a}{b}s^{2}}f(s)\right](y)j_{\nu}(\lambda y/b)y^{2\nu+1} dy\\
&=&e^{\frac{i}{2}(\frac{d}{b}\lambda^{2}-\frac{a}{b}x^{2})}
\ds\int_{0}^{+\infty}e^{\frac{i}{2}\frac{a}{b}y^{2}}f(y)T_{x}^{\nu }\left[s\mapsto j_{\nu}(\lambda s/b)\right](y)y^{2\nu+1} dy\\
&=&e^{\frac{i}{2}(\frac{d}{b}\lambda^{2}-\frac{a}{b}x^{2})}j_{\nu}(\lambda x/b)
\ds\int_{0}^{+\infty}e^{\frac{i}{2}\frac{a}{b}y^{2}}f(y)j_{\nu}(\lambda y/b)y^{2\nu+1} dy\\&=&
\left[(ib)^{\nu+1}/c_{\nu}\right] e^{\frac{i}{2}\frac{d}{b}\lambda^{2}}\ K_{\nu}^{\mathbf{m}^{-1}}(x,\lambda)\mathcal{F}_{\nu}^{\mathbf{m}}f(\lambda).
\end{eqnarray*}
$10)$ From $9)$ the result is true for $f\in \mathcal{L}_{1,\nu}\cap \mathcal{L}_{p,\nu}.$ On the other hand, Babenko inequality (\ref{babenko}) and the relation (\ref{inequa}) show that the mappings $f\mapsto\mathcal{F}_{\nu}^{\mathbf{m}}\left[T_{x}^{\nu,\mathbf{m}^{-1}}f\right]$ and
$f\mapsto\mathcal{F}_{\nu}^{\mathbf{m}}f$ are continuous from $\mathcal{L}_{p,\nu}$ into $\mathcal{L}_{q,\nu}\ (1/p+1/q=1).$
We obtain the result from density of $\mathcal{L}_{1,\nu}\cap \mathcal{L}_{p,\nu}$ in $\mathcal{L}_{p,\nu}.$\end{Proof-}
\begin{Corollary-} Let $\mathbf{m}\in SL(2,\mathbb{R})$ such that $b\neq0$ and $x\in\mathbb{R}.$ Then the operator $T_{x}^{\nu,\mathbf{m}^{-1}}$ leaves $\mathcal{S}_{\ast }(\mathbb{R})$ invariant and for each $f\in \mathcal{S}_{\ast }(\mathbb{R})$ we have:
\begin{eqnarray*}
T_{x}^{\nu,\mathbf{m}^{-1}}f(y)&=&\frac{c_{\nu}}{(-ib)^{\nu+1}}e^{-\frac{i}{2}\frac{a}{b}y^{2}}
\ds\int_{0}^{+\infty}j_{\nu}(\lambda y/b)K_{\nu}^{\mathbf{m}^{-1}}(x,\lambda)\mathcal{F}_{\nu}^{\mathbf{m}}f(\lambda)\lambda^{2\nu+1}\ d\lambda.
\end{eqnarray*}
\label{coro}
\end{Corollary-}
\begin{Proof-} To prove this, we first observe that if $f\in \mathcal{S}_{\ast }(\mathbb{R})$ then the inequality (\ref{inequa}) shows that
$y\mapsto[T_{x}^{\nu,\mathbf{m}^{-1}}f](y)$ is a continuous function of $ \mathcal{L}_{1,\nu}.$ The result is then a consequence of (\ref{fourier}) and the reversibility property of the canonical Fourier Bessel transform (see Theorem \ref{theorem1}).
\end{Proof-}

We conclude this subsection with the following two Theorems.
\begin{Theorem-} Let $\mathbf{m}\in SL(2,\mathbb{R})$ such that $b\neq0$.\\
$(i)$ For all $f\in \mathcal{C}_{*,0}(\mathbb{R}),$ we have $\ds\lim_{y\rightarrow 0}\left\|T_{y}^{\nu,\mathbf{m}}f-f\right\|_{\infty}=0.$\\
$(ii)$ For all $f\in \mathcal{L}_{p,\nu }$ with $1\leq p<+\infty,$ we have $\ds\lim_{y\rightarrow 0}\left\|T_{y}^{\nu,\mathbf{m}}f-f\right\|_{p,\nu}=0.$
\label{theo13}
\end{Theorem-}
\begin{Proof-} $(i)$ First consider the case where $f\in \mathcal{C}_{*,c}(\mathbb{R}).$ From the definition of $T_{y}^{\nu,\mathbf{m}}$ and the fact that $\frac{\Gamma(\nu+1)}{\sqrt{\pi}\Gamma(\nu+1/2)}\int_{0}^{\pi}\sin(\theta)^{2\nu} d\theta=1,$ we can write
\begin{eqnarray*}
T_{y}^{\nu,\mathbf{m}}f(x)-f(x)&=& T_{x}^{\nu,\mathbf{m}}f(y)-f(x)=a_{y}(x)+b_{y}(x)
\end{eqnarray*}
where
\begin{eqnarray*}
a_{y}(x)&=&f(x) \ \frac{\Gamma(\nu+1)}{\sqrt{\pi}\Gamma(\nu+1/2)}\int_{0}^{\pi}\left[e^{i\frac{d}{b}xy\cos(\theta)}-1\right]\sin(\theta)^{2\nu} d\theta,\\
b_{y}(x)&=&\frac{\Gamma(\nu+1)}{\sqrt{\pi}\Gamma(\nu+1/2)}\int_{0}^{\pi}e^{i\frac{d}{b}xy\cos(\theta)}
\left[f\left(\sqrt{x^{2}+y^{2}-2xy\cos(\theta)}\right)-f(x)\right]\sin(\theta)^{2\nu} d\theta.
\end{eqnarray*}
Clearly
\begin{eqnarray*}
a_{y}(x)&=&f(x) \ \frac{\Gamma(\nu+1)}{\sqrt{\pi}\Gamma(\nu+1/2)}\ds\sum_{n=1}^{+\infty}\frac{(i\frac{d}{b}xy)^{n}}{n!}\ds\int_{0}^{\pi}\cos^{n}(\theta)\left(\sin(\theta)\right)^{2\nu} d\theta,
\end{eqnarray*}
and since
\begin{eqnarray}
\ds\int_{0}^{\pi}\cos^{n}(\theta)\left(\sin(\theta)\right)^{2\nu} d\theta&=&[1+(-1)^{n}]\ds\int_{0}^{\frac{\pi}{2}}\cos^{n}(\theta)\left(\sin(\theta)\right)^{2\nu} d\theta=\left(\frac{1+(-1)^{n}}{2}\right)\frac{\Gamma(\frac{n+1}{2})\Gamma(\nu+\frac{1}{2})}{\Gamma(\frac{n}{2}+\nu +1)},
\label{gamma}
\end{eqnarray}
we conclude
\begin{eqnarray*}
a_{y}(x)&=&\left[j_{\nu}\left(\frac{dxy}{b}\right)-1\right]f(x).
\end{eqnarray*}
Let $R>0$ such that $\mbox{supp}f\subset [0,R].$ Then
$$\|a_{y}\|_{\infty}\leq \|f\|_{\infty}\left[j_{\nu}\left(i\frac{dRy}{b}\right)-1\right]\underset{y\rightarrow 0}{\longrightarrow 0}.
$$
By the triangle inequality, we have
\begin{eqnarray*}
\left|\sqrt{x^{2}+y^{2}-2xy\cos(\theta)}-x\right|&=&\left|\sqrt{(x-y\cos(\theta))^{2}+y^{2}\sin^{2}(\theta)}-\sqrt{x^{2}}\right|\leq |y|.
\end{eqnarray*}
Now given $\epsilon>0,$ by the uniform continuity of $f$ there exists $\delta>0$ such that $|f(u)-f(v)|\leq \epsilon$ whenever $|u-v|< \delta.$ Thus we get for $|y|< \delta$ and $x\in\mathbb{R},$
\begin{eqnarray*}
|b_{y}(x)|&\leq& \frac{\Gamma(\nu+1)}{\sqrt{\pi}\Gamma(\nu+1/2)}\int_{0}^{\pi}
\left|f(\sqrt{x^{2}+y^{2}-2xy\cos(\theta)})-f(x)\right|\sin(\theta)^{2\nu} d\theta \\
&\leq& \epsilon.
\end{eqnarray*}
Hence $$\ds\lim_{y\rightarrow 0}\left\|T_{y}^{\nu,\mathbf{m}}f-f\right\|_{\infty}=0.$$
Now suppose $f\in \mathcal{C}_{*,0}(\mathbb{R}).$ If $\epsilon>0$ and since $\mathcal{C}_{\ast,c}(\mathbb{R})$ is dense in $\mathcal{C}_{*,0}(\mathbb{R}),$ there exists $g\in \mathcal{C}_{\ast,c}(\mathbb{R})$ such that $\|f-g\|_{\infty}\leq \frac{\epsilon}{3},$ so
\begin{eqnarray*}
\left\|T_{y}^{\nu,\mathbf{m}}f-f\right\|_{\infty}&\leq&\left\|T_{y}^{\nu,\mathbf{m}}(f-g)\right\|_{\infty}
+\left\|T_{y}^{\nu,\mathbf{m}}g-g\right\|_{\infty}+\left\|g-f\right\|_{p,\nu}\\
&\leq&2\left\|g-f\right\|_{\infty}+\left\|T_{y}^{\nu,\mathbf{m}}g-g\right\|_{\infty}\leq \frac{2\epsilon}{3}+\left\|T_{y}^{\nu,\mathbf{m}}g-g\right\|_{\infty}
\end{eqnarray*}
and $\left\|T_{y}^{\nu,\mathbf{m}}g-g\right\|_{\infty}\leq \frac{\epsilon}{3}$ if $y$ is sufficiently small.\\
$(ii)$ First, if $f\in \mathcal{C}_{*,c}(\mathbb{R}),$ for $|y|\leq1$ the functions $T_{y}^{\nu,\mathbf{m}}f$ are all supported in a common compact set $[0,R],$ so
\begin{eqnarray*}
\left\|T_{y}^{\nu,\mathbf{m}}f-f\right\|_{p,\nu}^{p}\leq \left(\ds\int_{0}^{R}x^{2\nu+1}dx\right)\left\|T_{y}^{\nu,\mathbf{m}}f-f\right\|_{\infty}\rightarrow0 \ \mbox{as} \ y\rightarrow0.
\end{eqnarray*}
The general case then follows from density of $\mathcal{C}_{*,c}(\mathbb{R})$ in $\mathcal{L}_{p,\nu }.$
\end{Proof-}
\begin{Theorem-} Let $\mathbf{m}\in SL(2,\mathbb{R})$ such that $b\neq0.$ Let $x, y\in \mathbb{R}$ and $\sigma>0.$ For all $f\in \mathcal{S}_{\ast }(\mathbb{R}),$ we have:\\
$(i)$ The function $t\mapsto e^{-2i\frac{d}{b}\sigma ty^{2}}[T_{x}^{\nu,\mathbf{m}}f][\sqrt{4\sigma t}y]$ is $\mathcal{C}^{\infty}$ on $[0,+\infty[$ with
\begin{eqnarray*}
\frac{d^{n}}{dt^{n}}\left[ e^{-2i\frac{d}{b}\sigma ty^{2}}[T_{x}^{\nu,\mathbf{m}}f][\sqrt{4\sigma t}y]\right]
&=&\frac{\left(-\frac{\sigma}{b^{2}}\right)^{n}}{2^{\nu}(ib)^{\nu+1}\Gamma(n+\nu+1)}\\
& \times & y^{2n} \ds\int_{0}^{+\infty}j_{\nu+n}\left(\lambda \sqrt{4\sigma t}y/b\right)\lambda^{2n}K_{\nu}^{\mathbf{m}}(x,\lambda)\mathcal{F}_{\nu}^{\mathbf{m}^{-1}}f(\lambda)\lambda^{2\nu+1}\ d\lambda.
\end{eqnarray*}
$(ii)$ For all non-negative integers $n,$
\begin{eqnarray*}
\frac{d^{n}}{dt^{n}}\ \left.\left(e^{-2i\frac{d}{b}\sigma ty^{2}}[T_{x}^{\nu,\mathbf{m}}f][\sqrt{4\sigma t}y]\right)\right|_{t=0}=\frac{\Gamma(\nu+1)}{\Gamma(n+\nu+1)}\ \sigma^{n} \ y^{2n}\left(\Delta_{\nu}^{\mathbf{m}}\right)^{n}f(x).
\end{eqnarray*}
$(iii)$ For all non-negative integers $n$ and $t>0,$
\begin{eqnarray*}
e^{-2i\frac{d}{b}\sigma ty^{2}}[T_{x}^{\nu,\mathbf{m}}f][\sqrt{4\sigma t}y]&=&\ds\sum_{k=0}^{n}
\frac{\Gamma(\nu+1)}{\Gamma(k+\nu+1)k!}\ t^{k}\sigma^{k} \ y^{2k}\left(\Delta_{\nu}^{\mathbf{m}}\right)^{k}f(x)\\
&+&\frac{t^{n+1}}{n!}\ds\int_{0}^{1}(1-s)^{n}\frac{d^{n+1}}{du^{n+1}}\left[e^{-2i\frac{d}{b}\sigma uy^{2}}[T_{x}^{\nu,\mathbf{m}}f][\sqrt{4\sigma u}y]\right](st)\ ds.
\end{eqnarray*}
$(iv)$ For all non-negative integers $n,$ there is a constant $c_{n}>0$ and $n_{1},\dots ,n_{k_{n}}, m_{1},\dots, m_{k_{n}}\in \mathbb{N}$ such that
\begin{eqnarray*}
\left|\frac{d^{n}}{dt^{n}} e^{-2i\frac{d}{b}\sigma ty^{2}}[T_{x}^{\nu,\mathbf{m}}f][\sqrt{4\sigma t}y]\right|
&\leq& c_{n} \ y^{2n}\ds\sum_{i=1}^{k_{n}}q_{n_{i},m_{i}}(f).
\end{eqnarray*}
\label{deri}
\end{Theorem-}
\begin{Proof-}
$(i)$ By Corollary \ref{coro} we have
\begin{eqnarray*}
 e^{-2i\frac{d}{b}\sigma ty^{2}}[T_{x}^{\nu,\mathbf{m}}f][\sqrt{4\sigma t}y]&=&\frac{c_{\nu}}{(ib)^{\nu+1}}
 \ds\int_{0}^{+\infty}j_{\nu}\left(\lambda \sqrt{4\sigma t}y/b\right)K_{\nu}^{\mathbf{m}}(x,\lambda)\mathcal{F}_{\nu}^{\mathbf{m}^{-1}}f(\lambda)\lambda^{2\nu+1}\ d\lambda.
\end{eqnarray*}
The desired result is therefore a consequence of Theorem of differentiation under the integral sign; since
$$\frac{d^{n}}{dt^{n}}j_{\nu}\left(\lambda \sqrt{4\sigma t}y/b\right)=(-1)^{n}\frac{\Gamma(\nu+1)}{\Gamma(n+\nu+1)}\frac{\lambda^{2n} y^{2n}\sigma^{n}}{b^{2n}}j_{n+\nu}\left(\lambda \sqrt{4\sigma t}y/b\right),$$
$\left|j_{n+\nu}\left(\lambda \sqrt{4\sigma t}y/b\right)K_{\nu}^{\mathbf{m}}(x,\lambda)\right|\leq 1$ and $\lambda\mapsto \lambda^{2n}\mathcal{F}_{\nu}^{\mathbf{m}^{-1}}f\in \mathcal{S}_{\ast }(\mathbb{R}).$\\
$(ii)$ Clearly
\begin{eqnarray*}
\frac{d^{n}}{dt^{n}}\ \left.\left(e^{-2i\frac{d}{b}\sigma ty^{2}}[T_{x}^{\nu,\mathbf{m}}f][\sqrt{4\sigma t}y]\right)\right|_{t=0}&=&\frac{\left(-\frac{\sigma}{b^{2}}\right)^{n}y^{2n}}{2^{\nu}(ib)^{\nu+1}\Gamma(n+\nu+1)}\ds\int_{0}^{+\infty}
K_{\nu}^{\mathbf{m}}(x,\lambda)\lambda^{2n}\mathcal{F}_{\nu}^{\mathbf{m}^{-1}}f(\lambda)\lambda^{2\nu+1}\ d\lambda.
\end{eqnarray*}
By (\ref{equa1})
\begin{eqnarray*}
\lambda^{2n}\mathcal{F}_{\nu}^{\mathbf{m}^{-1}}f(\lambda)=
(-b^{2})^{n}\mathcal{F}_{\nu}^{\mathbf{m}^{-1}}\left[\left(\Delta_{\nu}^{\mathbf{m}}\right)^{n}f\right](\lambda).
\end{eqnarray*}
Then
\begin{eqnarray*}
\frac{d^{n}}{dt^{n}}\ \left.\left(e^{-2i\frac{d}{b}\sigma ty^{2}}[T_{x}^{\nu,\mathbf{m}}f][\sqrt{4\sigma t}y]\right)\right|_{t=0}&=&\frac{\Gamma(\nu+1)}{\Gamma(n+\nu+1)}\ \sigma^{n} \ y^{2n}\mathcal{F}_{\nu}^{\mathbf{m}}\left[\mathcal{F}_{\nu}^{\mathbf{m}^{-1}}\left(\Delta_{\nu}^{\mathbf{m}}\right)^{n}f\right]
(x)\\&=&\frac{\Gamma(\nu+1)}{\Gamma(n+\nu+1)}\ \sigma^{n} \ y^{2n}\left(\Delta_{\nu}^{\mathbf{m}}\right)^{n}f(x).
\end{eqnarray*}
$(iii)$ This follows from Taylor's formula and $(ii).$\\
$(iv)$ From $(i)$ it is clear that
\begin{eqnarray*}
\left|\frac{d^{n}}{dt^{n}} e^{-2i\frac{d}{b}\sigma ty^{2}}[T_{x}^{\nu,\mathbf{m}}f][\sqrt{4\sigma t}y]\right|
&\leq&\frac{\left(\frac{\sigma}{b^{2}}\right)^{n}}{2^{\nu}(|b|)^{\nu+1}\Gamma(n+\nu+1)}\\
& \times & y^{2n} \ds\int_{0}^{+\infty}\lambda^{2n}\left|\mathcal{F}_{\nu}^{\mathbf{m}^{-1}}f(\lambda)\right|\lambda^{2\nu+1}\ d\lambda.
\end{eqnarray*}
Now choose $m\in \mathbb{N}$ and $m>\nu+1.$ Then
\begin{eqnarray*}
\ds\int_{0}^{+\infty}\lambda^{2n}\left|\mathcal{F}_{\nu}^{\mathbf{m}^{-1}}f(\lambda)\right|\lambda^{2\nu+1}\ d\lambda&\leq& \left(\ds\int_{0}^{+\infty}\frac{\lambda^{2\nu+1}}{1+\lambda ^{2m}}d\lambda\right)
\sup_{\lambda\geq 0}\left((1+\lambda^{2m})\lambda^{2n}\left|\mathcal{F}_{\nu}^{\mathbf{m}^{-1}}f(\lambda) \right|\right)\\
&\leq&\left(\ds\int_{0}^{+\infty}\frac{\lambda^{2\nu+1}}{1+\lambda ^{2m}}d\lambda\right)\left(q_{n+m,0}[\mathcal{F}_{\nu}^{\mathbf{m}^{-1}}f]+q_{n,0}[\mathcal{F}_{\nu}^{\mathbf{m}^{-1}}f]\right).
\end{eqnarray*}
The desired result is therefore a consequence of (\ref{equat3}).\end{Proof-}

\subsection{Generalized convolution product}
Let $\mathbf{m}\in SL(2,\mathbb{R})$ such that $b\neq0$ and let $f$ and $g$ be measurable functions on $[0,+\infty[.$ The generalized convolution of $f$ and $g$ is the
function $f\underset{\nu,\mathbf{m}}{*}g$ defined by
\begin{eqnarray}
f\underset{\nu,\mathbf{m}}{*}g(x)=\int_{0}^{+\infty }\left[T_{x}^{\nu,\mathbf{m}}f
\right](y) \ \left[e^{-i\frac{d}{b}y^{2}}g(y)\right]y^{2\nu +1}dy
\label{convo}
\end{eqnarray}
for all $x$ such that the integral exists. The elementary properties of convolutions are summarized in the following
proposition.
\begin{Proposition-} Assuming that all integrals in question exists, we have\\
$1) \ f\underset{\nu,\mathbf{m}}{*}g=g\underset{\nu,\mathbf{m}}{*}f.$ \\
$2)\ T_{x}^{\nu,\mathbf{m}}\left(f\underset{\nu,\mathbf{m}}{*}g\right)=\left[T_{x}^{\nu,\mathbf{m}}f\right]\underset{\nu,\mathbf{m}}{*}g=f\underset{\nu,\mathbf{m}}{*}\left[T_{x}^{\nu,\mathbf{m}}g\right].$\\
$3) \left(f\underset{\nu,\mathbf{m}}{*}g\right)\underset{\nu,\mathbf{m}}{*}h=f\underset{\nu,\mathbf{m}}{*}\left(g\underset{\nu,\mathbf{m}}{*}h\right).$
\end{Proposition-}
\begin{Proof-}$1)$ This is (\ref{equat})

$2)$ By Fubini's theorem.

$3)$ follows from $(1), \ (2)$ and Fubini's theorem:
\end{Proof-}
\bigskip

The following proposition contain the basic facts about convolutions of $\mathcal{L}_{p,\nu}$
functions.
\begin{Proposition-} (\textbf{Young's Inequality}) Let $\mathbf{m}\in SL(2,\mathbb{R})$ such that $b\neq0.$ Suppose $1\leq p, q,r \leq\infty$ and $p^{-1}+q^{-1}=r^{-1}+1.$ If $f\in\mathcal{L}_{p,\nu}$ and $g\in\mathcal{L}_{q,\nu},$ then $f\underset{\nu,\mathbf{m}}{*}g \in \mathcal{L}_{r,\nu}$ and
\begin{eqnarray}
\left\|f\underset{\nu,\mathbf{m}}{*}g\right\|_{r,\nu}\leq \left\|f\right\|_{p,\nu} \ \left\|g\right\|_{q,\nu}.
\label{young}
\end{eqnarray}
\label{young1}
\end{Proposition-}
\begin{Proof-} By applying H\"{o}lder's inequality to the product
$$
\left|T_{x}^{\nu,\mathbf{m}}f(y)e^{-i\frac{d}{b}y^{2}}g(y)\right|=\left(\left|T_{x}^{\nu,\mathbf{m}}f(y)\right|^{p}|g(y)|^{q}\right)^{1/r}\left(\left|T_{x}^{\nu,\mathbf{m}}f(y)\right|^{p}\right)^{1/p-1/r}
\left(|g(y)|^{q}\right)^{1/q-1/r},
$$
we have
\begin{eqnarray*}
\int_{0}^{+\infty}\left|T_{x}^{\nu,\mathbf{m}}f(y)e^{-i\frac{d}{b}y^{2}}g(y)\right|y^{2\nu +1}dy&\leq&\left(\int_{0}^{+\infty}\left|T_{x}^{\nu,\mathbf{m}}f(y)\right|^{p}|g(y)|^{q}y^{2\nu +1}dy\right)^{\frac{1}{r}}
\left(\int_{0}^{+\infty}\left|T_{x}^{\nu,\mathbf{m}}f(y)\right|^{p}y^{2\nu +1}dy\right)^{\frac{r-p}{rp}}\\
&\times&\left(\int_{0}^{+\infty}|g(y)|^{q}y^{2\nu +1}dy\right)^{\frac{r-q}{rq}}.
\end{eqnarray*}
This leads to:
\begin{eqnarray*}
\left|\left(f\underset{\nu,\mathbf{m}}{*}g\right)(x)\right|^{r}&\leq&\left\|T_{x}^{\nu,\mathbf{m}}f\right\|_{p,\nu} ^{r-p} \left\|g\right\|_{q,\nu} ^{r-q} \int_{0}^{+\infty}\left|T_{x}^{\nu,\mathbf{m}}f(y)\right|^{p}|g(y)|^{q}y^{2\nu +1}dy\\
\end{eqnarray*}
and using (\ref{inequa}) we obtain
\begin{eqnarray*}
\left|\left(f\underset{\nu,\mathbf{m}}{*}g\right)(x)\right|^{r}&\leq& \left\|f\right\|_{p,\nu} ^{r-p} \left\|g\right\|_{q,\nu} ^{r-q} \int_{0}^{+\infty}\left|T_{x}^{\nu,\mathbf{m}}f(y)\right|^{p}|g(y)|^{q}y^{2\nu +1}dy.
\end{eqnarray*}
Multiply both sides by $x^{2\nu+1}$ and integrate from $0$ to $+\infty,$ we get
\begin{eqnarray*}
\left\|f\underset{\nu,\mathbf{m}}{*}g\right\|_{r,\nu}^{r}&\leq&  \left\|f\right\|_{p,\nu} ^{r-p} \left\|g\right\|_{q,\nu} ^{r-q} \int_{0}^{+\infty} \left[\int_{0}^{+\infty}\left|T_{x}^{\nu,\mathbf{m}}f(y)\right|^{p}|g(y)|^{q}y^{2\nu +1}dy\right]x^{2\nu+1} dx\\
&=&\left\|f\right\|_{p,\nu} ^{r-p} \left\|g\right\|_{q,\nu} ^{r-q} \int_{0}^{+\infty} |g(y)|^{q}y^{2\nu +1}\left[\int_{0}^{+\infty}\left|T_{y}^{\nu,\mathbf{m}}f(x)\right|^{p}x^{2\nu+1} dx\right]\\
&\leq&\left\|f\right\|_{p,\nu} ^{r} \left\|g\right\|_{q,\nu} ^{r}.
\end{eqnarray*}
\end{Proof-}
\begin{Proposition-} Let $\mathbf{m}\in SL(2,\mathbb{R})$ such that $b\neq0.$\\
$1)$ Let $f$ and $g$ be two function in $\mathcal{L}_{1,\nu}.$ We have
\begin{eqnarray*}
\forall x\in \mathbb{R}, \quad \left(c_{\nu}/(ib)^{\nu+1}\right)\mathcal{F}_{\nu}^{\mathbf{m}}\left(f\underset{\nu,\mathbf{m}^{-1}}{*}g\right)(x)=e^{-\frac{i}{2}\frac{d}{b}x^{2}} \  \mathcal{F}_{\nu}^{\mathbf{m}}f(x) \ \mathcal{F}_{\nu}^{\mathbf{m}}g(x).
\end{eqnarray*}
$2)$ Let $f\in \mathcal{L}_{1,\nu} $ and $g\in\mathcal{L}_{p,\nu}\ (p\in]1,2]).$ We have
\begin{eqnarray*}
\left(c_{\nu}/(ib)^{\nu+1}\right)\mathcal{F}_{\nu}^{\mathbf{m}}\left(f\underset{\nu,\mathbf{m}^{-1}}{*}g\right)(x)=e^{-\frac{i}{2}\frac{d}{b}x^{2}} \  \mathcal{F}_{\nu}^{\mathbf{m}}f(x) \ \mathcal{F}_{\nu}^{\mathbf{m}}g(x), \ \mbox{a.e.}
\end{eqnarray*}
\label{propo}
\end{Proposition-}
\begin{Proof-}$1)$ According to the previous proposition, one has $(f\underset{\nu,\mathbf{m}^{-1}}{*}g)\in \mathcal{L}_{1,\nu}.$ From the definition of $\mathcal{F}_{\nu}^{\mathbf{m}}$ and (\ref{fourier}) it is easy to see that
\begin{eqnarray*}
 \mathcal{F}_{\nu}^{\mathbf{m}}\left(f\underset{\nu,\mathbf{m}^{-1}}{*}g\right)(x)&=&\frac{c_{\nu}}{(ib)^{\nu+1}}
 \int_{0}^{+\infty}K_{\nu}^{\mathbf{m}}(x,y)\left[\int_{0}^{+\infty}T_{y}^{\nu,\mathbf{m}^{-1}}f(z)\left[e^{i\frac{a}{b}z^{2}}g(z)\right]z^{2\nu+1}dz \right]y^{2\nu+1}dy\\
 &=&\frac{c_{\nu}}{(ib)^{\nu+1}}
 \int_{0}^{+\infty}\left[e^{i\frac{a}{b}z^{2}}g(z)\right]\left[\int_{0}^{+\infty}K_{\nu}^{\mathbf{m}}(x,y)T_{z}^{\nu,\mathbf{m}^{-1}}f(y)y^{2\nu+1}dy \right]z^{2\nu+1}dz\\&=&
 \int_{0}^{+\infty}\left[e^{i\frac{a}{b}z^{2}}g(z)\right]\left[ \mathcal{F}_{\nu}^{\mathbf{m}}\left(T_{z}^{\nu,\mathbf{m}^{-1}}f\right)(x)\right]z^{2\nu+1}dz\\
 &=&\left((ib)^{\nu+1}/c_{\nu}\right)\ e^{-\frac{i}{2}\frac{d}{b}x^{2}}\mathcal{F}_{\nu}^{\mathbf{m}}f(x)\mathcal{F}_{\nu}^{\mathbf{m}}g(x),
\end{eqnarray*}
The Fubini theorem was used in the second line, since
\begin{eqnarray*}
\int_{0}^{+\infty}\int_{0}^{+\infty}\left|K_{\nu}^{\mathbf{m}}(x,y)T_{y}^{\nu,\mathbf{m}^{-1}}f(z)e^{i\frac{a}{b}z^{2}}g(z)\right|y^{2\nu+1}z^{2\nu+1}dydz
&\leq&\int_{0}^{+\infty}\int_{0}^{+\infty}\left|T_{y}^{\nu,\mathbf{m}^{-1}}f(z)\right|\ |g(z)|y^{2\nu+1}z^{2\nu+1}dydz\\
&=&\|f\|_{1,\nu} \ \|g\|_{1,\nu}<\infty.
\end{eqnarray*}
$2)$ From $1)$ the result is true for $g\in \mathcal{L}_{1,\nu}\cap \mathcal{L}_{p,\nu}.$ On the other hand, the Babenko inequality (\ref{babenko}) and Proposition \ref{young1} show that the mappings $g\mapsto\mathcal{F}_{\nu}^{\mathbf{m}}\left(f\underset{\nu,\mathbf{m}^{-1}}{*}g\right)$ and
$g\mapsto\mathcal{F}_{\nu}^{\mathbf{m}}f \ \mathcal{F}_{\nu}^{\mathbf{m}}g$ are continuous from $\mathcal{L}_{p,\nu}$ into $\mathcal{L}_{q,\nu}\ (1/p+1/q=1).$
We obtain the result from density of $\mathcal{L}_{1,\nu}\cap \mathcal{L}_{p,\nu}$ in $\mathcal{L}_{p,\nu}.$
\end{Proof-}
\section{Heat equations and heat semigroups related to $\Delta_{\nu}^{\mathbf{m}^{-1}}$.}
\paragraph{}
In this section, we shall illustrate our theory by studying the following heat equations associated to the operator $\Delta_{\nu}^{\mathbf{m}^{-1}}:$
\begin{eqnarray}
\left\{\begin{array}{l}
\frac{\partial}{\partial t}u(t,x)=\sigma \Delta_{\nu}^{\mathbf{m}^{-1}}u(t,x), \ (t,x)\in (0,+\infty)\times\mathbb{R}  \\
u(0,x)=f(x),
\end{array}
\right.
\label{cauchy}
\end{eqnarray}
where $f$ is a given function in $\mathcal{L}_{p,\nu} \ (1\leq p\leq \infty)$ and $\sigma>0$ is the coefficient of heat conductivity. Here the initial data $u(0,x)=f(x)$ means that $u(t,x)\rightarrow f(x)$ as $t\rightarrow 0$ in the norm of $\mathcal{L}_{p,\nu}$.
\paragraph{}
In order to solve the Cauchy problem (\ref{cauchy}), it suggests itself to apply the canonical Fourier Bessel transform $\mathcal{F}_{\nu}^{\mathbf{m}}$, thus converting the partial differential equation $(\partial_{t}-\sigma \Delta_{\nu}^{\mathbf{m}^{-1}})u=0$ into the simple ordinary differential equation $[\partial_{t}+(\sigma x^{2}/b^{2})]\mathcal{F}_{\nu}^{\mathbf{m}}\left[s\mapsto u(t,s)\right](x)=0.$ The general solution of this equation is
\begin{eqnarray*}
\mathcal{F}_{\nu}^{\mathbf{m}}\left[s\mapsto u(t,s)\right](x)=c(x)e^{-\frac{\sigma t}{b^{2}}x^{2}},
\end{eqnarray*}
and we require that $\mathcal{F}_{\nu}^{\mathbf{m}}\left[s\mapsto u(0,s)\right](x)=\mathcal{F}_{\nu}^{\mathbf{m}}f(x).$ We therefore obtain a solution to our problem
by taking $c(x)=\mathcal{F}_{\nu}^{\mathbf{m}}f(x);$ this gives
$
\mathcal{F}_{\nu}^{\mathbf{m}}\left[s\mapsto u(t,s)\right](x)=e^{-\frac{\sigma t}{b^{2}}x^{2}}\mathcal{F}_{\nu}^{\mathbf{m}}f(x).
$
By
\begin{eqnarray*}
\ds\int_{0}^{+\infty}e^{-\beta y^{2}}j_{\nu}\left(\frac{xy}{b}\right)y^{2\nu+1} \ dy&=&
\frac{\Gamma(\nu+1)}{2\beta^{\nu+1}}e^{-\frac{x^{2}}{4\beta b^{2}}}; \ \beta>0, \ b\neq 0,
\end{eqnarray*}
we can write
\begin{eqnarray}
e^{-\frac{\sigma t}{b^{2}}x^{2}}&=&\left[(ib)^{\nu+1}/c_{\nu}\right]e^{-i\frac{d}{b}x^{2}}\mathcal{F}_{\nu}^{\mathbf{m}}\left[\mathcal{P}_{t}^{\mathbf{m}^{-1}}\right](x),
\label{the formula}
\end{eqnarray}
where
$$\mathcal{P}_{t}^{\mathbf{m}^{-1}}(x)=\left[2/\Gamma(\nu+1)\right](4\sigma t)^{-(\nu+1)}e^{-\frac{i}{2}\frac{a}{b}x^{2}-\frac{x^{2}}{4\sigma t}}.$$
After applying the Proposition \ref{propo}, we obtain
\begin{eqnarray*}
\mathcal{F}_{\nu}^{\mathbf{m}}\left[u(t,s)\right](x)&=&\left[(ib)^{\nu+1}/c_{\nu}\right]e^{-i\frac{d}{b}x^{2}}\mathcal{F}_{\nu}^{\mathbf{m}}
\left[\mathcal{P}_{t}^{\mathbf{m}^{-1}}\right](x)\mathcal{F}_{\nu}^{\mathbf{m}}f(x)\\
&=&\mathcal{F}_{\nu}^{\mathbf{m}}[\mathcal{P}_{t}^{\mathbf{m}^{-1}}\underset{\nu,\mathbf{m}^{-1}}{*}f](x).
\end{eqnarray*}
So we have a candidate for a solution:
\begin{eqnarray*}
u(t,x)=[\mathcal{P}_{t}^{\mathbf{m}^{-1}}\underset{\nu,\mathbf{m}^{-1}}{*}f](x)=\int_{0}^{+\infty }T_{x}^{\nu,\mathbf{m}^{-1}}[\mathcal{P}_{t}^{\mathbf{m}^{-1}}](y) \ [e^{i\frac{a}{b}y^{2}}f(y)]\ y^{2\nu +1}dy.
\end{eqnarray*}
So far this is all formal, since we have not specified conditions on $f$ to ensure that these manipulations are justified. Our object
is to show that the function $u(t,x)=[\mathcal{P}_{t}^{\mathbf{m}^{-1}}\underset{\nu,\mathbf{m}^{-1}}{*}f](x)$
with $(t,x)\in (0,+\infty)\times\mathbb{R}$ and $f\in\mathcal{L}_{p,\nu}$ solves the heat equation $(\partial_{t}-\sigma \Delta_{\nu}^{\mathbf{m}^{-1}})u=0$ and we study the problem of uniqueness of solutions of the Cauchy problem (\ref{cauchy}).\\

We begin with the following Lemma.
\begin{Lemma-} Let $\delta , r , s \in\mathbb{C}$ such that $\mbox{Rel}(\delta)>0.$ Then
\begin{eqnarray*}
\ds\int_{0}^{+\infty}e^{-\delta x^{2}}j_{\nu}(2rx)j_{\nu}(2sx) x^{2\nu+1}\ dy=\frac{\Gamma(\nu+1)}{2\delta^{\nu+1}} \ e^{-(1/\delta)(s^{2}+r^{2})}j_{\nu}(2irs/\delta).
\end{eqnarray*}
\label{lemma1}
\end{Lemma-}
\begin{Proof-}
We need the following formulas (see 6.615 in \cite{Gradshteyn}):
\begin{eqnarray}
\ds\int_{0}^{+\infty}e^{-\delta y}J_{\nu}(2r\sqrt{y})J_{\nu}(2s\sqrt{y})\ dy=\frac{e^{-(1/\delta)(s^{2}+r^{2})}}{\delta}I_{\nu}(2rs/\delta); \quad  |\arg(\delta)|< \pi/2, \quad \nu\geq 0,\label{equation24}
\end{eqnarray}
where
\begin{eqnarray}
I_{\nu}(z)=\ds\sum_{n=0}^{+\infty}\frac{(z/2)^{2n+\nu}}{n! \ \Gamma(n+\nu+1)}. \label{equation55}
\end{eqnarray}
In view of the relationship
\begin{eqnarray*}
j_{\mu}(x)=2^{\mu}\Gamma(\mu+1)\ds\frac{J_{\mu}(x)}{x^{\mu}}
\end{eqnarray*}
between the normalized Bessel function $j_{\mu}$ and the classical Bessel function $J_{\mu},$ the equation (\ref{equation24}) takes when $\delta, r, s>0$ the following form:
\begin{eqnarray}
\ds\int_{0}^{+\infty}e^{-\delta y}j_{\nu}(2r\sqrt{y})j_{\nu}(2s\sqrt{y}) y^{\nu}\ dy=\frac{\Gamma(\nu+1)}{\delta^{\nu+1}} \ e^{-(1/\delta)(s^{2}+r^{2})}j_{\nu}(2irs/\delta).\label{equation25}
\end{eqnarray}
After the change of variables $x=\sqrt{y}$ and the analytic continuation of holomorphic functions, the equation (\ref{equation25}) becomes
\begin{eqnarray*}
\ds\int_{0}^{+\infty}e^{-\delta x^{2}}j_{\nu}(2rx)j_{\nu}(2sx) x^{2\nu+1}\ dy=\frac{\Gamma(\nu+1)}{2\delta^{\nu+1}} \ e^{-(1/\delta)(s^{2}+r^{2})}j_{\nu}(2irs/\delta); \quad  \delta , r , s \in\mathbb{C} , \ \mbox{Rel}(\delta)>0.
\end{eqnarray*}
\end{Proof-}
\begin{Definition-} Let $\mathbf{m}\in SL(2,\mathbb{R})$ such that $b\neq0.$ The generalized heat kernel $G_{t}^{\mathbf{m}^{-1}}$ associated to $\Delta_{\nu}^{\mathbf{m}^{-1}}$ is defined by
\begin{eqnarray}
G_{t}^{\mathbf{m}^{-1}}(x,y)=T_{x}^{\nu,\mathbf{m}^{-1}}[\mathcal{P}_{t}^{\mathbf{m}^{-1}}](y), \ x,y\in \mathbb{R}, \ t>0.
\label{heat}
\end{eqnarray}
\end{Definition-}
We collect some basic properties of the generalized heat kernel $G_{t}^{\mathbf{m}^{-1}}.$
\begin{Proposition-} Let $\mathbf{m}\in SL(2,\mathbb{R})$ such that $b\neq0$. The generalized heat kernel $G_{t}^{\mathbf{m}^{-1}}$ has the following properties:\\
$1)$ \ $G_{t}^{\mathbf{m}^{-1}}(x,y)=\frac{2}{\Gamma(\nu+1)(4\sigma t)^{\nu+1}}e^{-\frac{i}{2}\frac{a}{b}(x^{2}+y^{2})}e^{-\frac{x^{2}+y^{2}}{4\sigma t}}j_{\nu}\left(\frac{i xy}{2\sigma t}\right).$\\
$2)$ \ $|G_{t}^{\mathbf{m}^{-1}}(x,y)|\leq \frac{2 e^{-\frac{(|x|-|y|)^{2}}{4\sigma t}}}{\Gamma(\nu+1)(4\sigma t)^{\nu+1}}.$\\
$3)$ \ $\ds\int_{0}^{+\infty }e^{\frac{i}{2}\frac{a}{b}(x^{2}+y^{2})}G_{t}^{\mathbf{m}^{-1}}(x,y)y^{2\nu +1}dy=1.$\\
$4)$ \ $G_{t+s}^{\mathbf{m}^{-1}}(x,y)=\ds\int_{0}^{+\infty }G_{t}^{\mathbf{m}^{-1}}(x,z)G_{s}^{\mathbf{m}^{-1}}(y,z)\ e^{i\frac{a}{b}z^{2}}z^{2\nu +1}dz.$\\
$5)$ For fixed $y\in \mathbb{R},$ the function $u(t,x)=G_{t}^{\mathbf{m}^{-1}}(x,y)$ solves the heat equation $\frac{\partial}{\partial t}u(t,x)=\sigma \Delta_{\nu}^{\mathbf{m}^{-1}}u(t,x)$ on $(0,+\infty)\times\mathbb{R}.$ ($G_{t}^{\mathbf{m}^{-1}}$ defines a fundamental solution of the heat equation (\ref{cauchy})).
\label{semigroup}
\end{Proposition-}
\begin{Proof-}$1)$ From the definition of $T_{x}^{\nu,\mathbf{m}^{-1}}$ it follows that
\begin{eqnarray*}
G_{t}^{\mathbf{m}^{-1}}(x,y)&=&\frac{2}{\Gamma(\nu+1)(4\sigma t)^{\nu+1}}e^{-\frac{i}{2}\frac{a}{b}(x^{2}+y^{2})}
T_{x}^{\nu}[e^{-\frac{s^{2}}{4\sigma t}}](y).
\end{eqnarray*}
Now a simple calculation shows that
\begin{eqnarray*}
T_{x}^{\nu}[e^{-\frac{s^{2}}{4\sigma t}}](y)&=&\frac{\Gamma(\nu+1)}{\sqrt{\pi}\Gamma(\nu+1/2)}e^{-\frac{x^{2}+y^{2}}{4\sigma t}}\ds\int_{0}^{\pi}e^{\frac{xy}{2\sigma t}\cos(\theta)}\left(\sin(\theta)\right)^{2\nu} d\theta\\
&=&\frac{\Gamma(\nu+1)}{\sqrt{\pi}\Gamma(\nu+1/2)}e^{-\frac{x^{2}+y^{2}}{4\sigma t}}\ds\sum_{n=0}^{+\infty}
\frac{(\frac{xy}{2\sigma t})^{n}}{n!}\ds\int_{0}^{\pi}\cos^{n}(\theta)\left(\sin(\theta)\right)^{2\nu} d\theta.
\end{eqnarray*}
By (\ref{gamma}) we conclude
\begin{eqnarray*}
T_{x}^{\nu}[e^{-\frac{s^{2}}{4\sigma t}}](y)=e^{-\frac{x^{2}+y^{2}}{4\sigma t}}j_{\nu}\left(\frac{i xy}{2\sigma t}\right).
\end{eqnarray*}
$3)$ From Lemma \ref{lemma1} we know that
\begin{eqnarray*}
\ds\int_{0}^{+\infty }e^{\frac{i}{2}\frac{a}{b}(x^{2}+y^{2})}G_{t}^{\mathbf{m}^{-1}}(x,y)y^{2\nu +1}dy
&=&\frac{2 e^{-\frac{x^{2}}{4\sigma t}}}{\Gamma(\nu+1)(4\sigma t)^{\nu+1}}
\ds\int_{0}^{+\infty }e^{-\frac{y^{2}}{4\sigma t}}j_{\nu}\left(\frac{i xy}{2\sigma t}\right)y^{2\nu +1}dy\\
&=&1.
\end{eqnarray*}
$4)$ We have
\begin{eqnarray*}
\ds\int_{0}^{+\infty }G_{t}^{\mathbf{m}^{-1}}(x,z)G_{s}^{\mathbf{m}^{-1}}(y,z)\ e^{i\frac{a}{b}z^{2}}z^{2\nu +1}dz&=&\frac{4 e^{-\frac{i}{2}\frac{a}{b}(x^{2}+y^{2})}e^{-\frac{x^{2}}{4\sigma t}}e^{-\frac{y^{2}}{4\sigma s}}}{\Gamma^{2}(\nu+1)(4\sigma t)^{\nu+1}(4\sigma s)^{\nu+1}}\\
&\times&\ds\int_{0}^{+\infty }e^{-\frac{z^{2}}{4\sigma }\frac{t+s}{ts}}j_{\nu}\left(\frac{i xz}{2\sigma t}\right)j_{\nu}\left(\frac{i yz}{2\sigma s}\right)z^{2\nu +1}dz.
\end{eqnarray*}
From Lemma \ref{lemma1}
\begin{eqnarray*}
\ds\int_{0}^{+\infty }e^{-\frac{z^{2}}{4\sigma }\frac{t+s}{ts}}j_{\nu}\left(\frac{i xz}{2\sigma t}\right)j_{\nu}\left(\frac{i yz}{2\sigma s}\right)z^{2\nu +1}dz=\frac{(4\sigma ts)^{\nu+1}\Gamma(\nu+1)}{2(t+s)^{\nu+1}}e^{\frac{sx^{2}}{4\sigma t(t+s)}+\frac{ty^{2}}{4\sigma s(t+s)}}
j_{\nu}\left(\frac{i xy}{2\sigma (t+s)}\right).
\end{eqnarray*}
These give the desired result.\\
$5)$ We simply compute for $y\in\mathbb{R}$ and $t>0,$
\begin{eqnarray}
\frac{\partial}{\partial t}G_{t}^{\mathbf{m}^{-1}}(x,y)&=&\left[\frac{x^{2}+y^{2}}{4\sigma t^{2}}-\frac{\nu+1}{t}\right]G_{t}^{\mathbf{m}^{-1}}(x,y)+\frac{2e^{-\frac{i}{2}\frac{a}{b}(x^{2}+y^{2})}e^{-\frac{x^{2}+y^{2}}{4\sigma t}}}{\Gamma(\nu+1)(4\sigma t)^{\nu+1}}\left(-\frac{i xy}{2\sigma t^{2}}\right)j_{\nu}^{'}\left(\frac{i xy}{2\sigma t}\right).
\label{partial}
\end{eqnarray}
By the transmutation property
\begin{eqnarray*}
e^{-\frac{i}{2}\frac{d}{b}x^{2}}\circ \Delta_{\nu}^{\mathbf{m}}\ \circ e^{\frac{i}{2}\frac{d}{b}x^{2}}=\Delta_{\nu},
\end{eqnarray*}
we have
\begin{eqnarray*}
\sigma\Delta_{\nu}^{\mathbf{m}^{-1}}\left[G_{t}^{\mathbf{m}^{-1}}(x,y)\right]&=&\sigma e^{-\frac{i}{2}\frac{a}{b}x^{2}}\Delta_{\nu}\left[e^{\frac{i}{2}\frac{a}{b}x^{2}}G_{t}^{\mathbf{m}^{-1}}(x,y)\right]\\
&=&\frac{2\sigma e^{-\frac{i}{2}\frac{a}{b}(x^{2}+y^{2})}e^{-\frac{y^{2}}{4\sigma t}}}{\Gamma(\nu+1)(4\sigma t)^{\nu+1}}\Delta_{\nu}\left[e^{-\frac{x^{2}}{4\sigma t}}j_{\nu}\left(\frac{i xy}{2\sigma t}\right)\right].
\end{eqnarray*}
An easy calculation shows that
\begin{eqnarray*}
\Delta_{\nu}\left[e^{-\frac{x^{2}}{4\sigma t}}j_{\nu}\left(\frac{i xy}{2\sigma t}\right)\right]&=&e^{-\frac{x^{2}}{4\sigma t}}\left[\left(\frac{iy}{2\sigma t}\right)^{2}j_{\nu}^{''}\left(\frac{i xy}{2\sigma t}\right)+\left(\frac{(2\nu+1)iy}{2\sigma t x}-\frac{ixy}{2\sigma^{2} t^{2}}\right)j_{\nu}^{'}\left(\frac{i xy}{2\sigma t}\right)\right.\\
&=&\left.+\left(\frac{x^{2}}{4\sigma^{2}t^{2}}-\frac{\nu+1}{\sigma t}\right)j_{\nu}\left(\frac{i xy}{2\sigma t}\right)\right].
\end{eqnarray*}
Since $j_{\nu}^{''}(z)+\frac{2\nu+1}{z}j_{\nu}^{'}(z)=-j_{\nu}(z), \ z\in \mathbb{C},$ we deduce
\begin{eqnarray*}
\Delta_{\nu}\left[e^{-\frac{x^{2}}{4\sigma t}}j_{\nu}\left(\frac{i xy}{2\sigma t}\right)\right]&=&e^{-\frac{x^{2}}{4\sigma t}}\left[\left(\frac{x^{2}+y^{2}}{4\sigma^{2} t^{2}}-\frac{\nu+1}{t}\right)j_{\nu}\left(\frac{i xy}{2\sigma t}\right)-\frac{ixy}{2\sigma^{2} t^{2}}j_{\nu}^{'}\left(\frac{i xy}{2\sigma t}\right)\right].
\end{eqnarray*}
Then
\begin{eqnarray*}
\left[\frac{\partial}{\partial t}-\sigma\Delta_{\nu}^{\mathbf{m}^{-1}}\right]G_{t}^{\mathbf{m}^{-1}}(x,y)=0.
\end{eqnarray*}\end{Proof-}
\begin{Definition-} Let $\mathbf{m}\in SL(2,\mathbb{R})$ such that $b\neq0.$ For $f\in \mathcal{L}_{p,\nu}$ and $t\geq0,$ set
\begin{eqnarray*}
S_{\nu}^{\mathbf{m}^{-1}}(t)f=\left\{\begin{array}{ll}
[\mathcal{P}_{t}^{\mathbf{m}^{-1}}\underset{\nu,\mathbf{m}^{-1}}{*}f] & \ \mbox{if} \ t>0, \\
f & \ \mbox{if} \ t=0.
\end{array}
\right.
\end{eqnarray*}
\end{Definition-}
\begin{Theorem-}
For each $f\in\mathcal{L}_{p,\nu} \ (1\leq p\leq \infty)$ set $u(t,x)=S_{\nu}^{\mathbf{m}^{-1}}(t)f(x)$ for any $(t,x)\in(0,+\infty)\times\mathbb{R}.$ Then:\\
$(1)$ The function $u$ satisfies the heat equation $(\partial_{t}-\sigma \Delta_{\nu}^{\mathbf{m}^{-1}})u=0$ in $(0,+\infty)\times\mathbb{R}.$\\
$(2)$ The solution $u(t,x)\rightarrow f(x)$ as $t\rightarrow 0$ in the norm of $\mathcal{L}_{p,\nu}$.\\
$(3)$\ $ \|u(t,\cdot)\|_{r,\nu }\leq \left[(4\sigma t)^{\nu+1}\Gamma(\nu+1)/2\right]^{1/q-1}q^{-(\nu+1)/q} \ \|f\|_{p,\nu }$
with $p,q,r\in[1,+\infty]$ be such that $1/p+1/q=1+1/r$ and $t>0.$\\
$(4)$ $u(t,.)$ is the unique solution of (\ref{cauchy}) belonging in $\mathcal{L}_{r,\nu}$ with $r\in[1,+\infty]$.\\
$(5)$ If $f\in \mathcal{S}_{\ast }(\mathbb{R})$ then $u(t,.)\in \mathcal{S}_{\ast }(\mathbb{R})$ and for all non-negative integers $n,$ there is a constant $c_{n+1}>0$ and $n_{1},\dots ,n_{k_{n+1}}$, $m_{1},\dots, m_{k_{n+1}}\in \mathbb{N}$ such that
\begin{eqnarray*}
\left\|u(t,.)-\ds\sum_{k=0}^{n}\frac{t^{k}}{k!} \left(\sigma\Delta_{\nu}^{\mathbf{m}^{-1}}\right)^{k}f\right\|_{\infty}\leq \frac{t^{n+1}}{(n+1)!}c_{n+1} \ds\sum_{i=1}^{k_{n+1}}q_{n_{i},m_{i}}(f).
\end{eqnarray*}
\label{proposi}
\end{Theorem-}
\begin{Proof-} $1)$ Let $f\in\mathcal{L}_{p,\nu} \ (1\leq p\leq \infty)$ In view of (\ref{convo}) and (\ref{heat}) we can write
$$u(t,x)=[\mathcal{P}_{t}^{\mathbf{m}^{-1}}\underset{\nu,\mathbf{m}^{-1}}{*}f](x)=\ds\int_{0}^{+\infty}G_{t}^{\mathbf{m}^{-1}}(x,y)
e^{i\frac{a}{b}y^{2}}f(y)y^{2\nu+1}\ dy.
$$
For any $(t,x)\in]0,+\infty[\times \mathbb{R}$ the function $y\mapsto G_{t}^{\mathbf{m}^{-1}}(x,y)$ is in $\mathcal{L}_{q,\nu}$ for $1\leq q\leq \infty.$ Then $[\mathcal{P}_{t}^{\mathbf{m}^{-1}}\underset{\nu,\mathbf{m}^{-1}}{*}f]$ is well defined and by (\ref{partial}), we can show that
\begin{eqnarray*}
\frac{\partial}{\partial t}u(t,x)= \ds\int_{0}^{+\infty}\frac{\partial}{\partial t}G_{t}^{\mathbf{m}^{-1}}(x,y)
e^{i\frac{a}{b}y^{2}}f(y)y^{2\nu+1}\ dy.
\end{eqnarray*}
Now a rough calculation shows that
\begin{eqnarray*}
\left|\sigma \Delta_{\nu}^{\mathbf{m}^{-1}} G_{t}^{\mathbf{m}^{-1}}(x,y)\right|&\leq& \frac{2}{\Gamma(\nu+1)(4\sigma t)^{\nu+1}}\left(\frac{r^{2}+y^{2}}{4\sigma t^{2}}+\frac{\nu+1}{t}+\frac{ry}{2\sigma t^{2}}\right)e^{-\frac{(y-r)^{2}}{4\sigma t}}
\end{eqnarray*}
for $|x|\leq r$ and $y\geq r.$ Hence
\begin{eqnarray*}
\sigma \Delta_{\nu}^{\mathbf{m}^{-1}}u(t,x)= \ds\int_{0}^{+\infty}\sigma \Delta_{\nu}^{\mathbf{m}^{-1}}G_{t}^{\mathbf{m}^{-1}}(x,y) \
e^{i\frac{a}{b}y^{2}}f(y)y^{2\nu+1}\ dy
\end{eqnarray*}
which yields the desired result by virtue of $\frac{\partial}{\partial t}G_{t}^{\mathbf{m}^{-1}}(x,y)=\sigma \Delta_{\nu}^{\mathbf{m}^{-1}}G_{t}^{\mathbf{m}^{-1}}(x,y).$\\
$2)$ We begin by showing that if $f\in \mathcal{C}_{\ast,c}(\mathbb{R}),$ then
$$\ds\lim_{t\rightarrow 0}\left\|S_{\nu}^{\mathbf{m}^{-1}}(t)f-f\right\|_{p,\nu}=0.$$
By (\ref{equat}), we can write
\begin{eqnarray}
S_{\nu}^{\mathbf{m}^{-1}}(t)f(x)&=&\int_{0}^{+\infty}[T_{x}^{\nu,\mathbf{m}^{-1}}\mathcal{P}_{t}^{\mathbf{m}^{-1}}](y)[e^{i\frac{a}{b}y^{2}}
f(y)]y^{2\nu+1}dy\nonumber\\&=&\int_{0}^{+\infty}e^{i\frac{a}{b}y^{2}}\mathcal{P}_{t}^{\mathbf{m}^{-1}}(y)\ [T_{x}^{\nu,\mathbf{m}^{-1}}f](y)y^{2\nu+1}dy\nonumber\\&=&\frac{2}{\Gamma(\nu+1)(4\sigma t)^{\nu+1}}\int_{0}^{+\infty}e^{\frac{i}{2}\frac{a}{b}y^{2}}e^{-\frac{y^{2}}{4\sigma t}}\ [T_{x}^{\nu,\mathbf{m}^{-1}}f](y)y^{2\nu+1}dy\nonumber\\&=&\frac{2}{\Gamma(\nu+1)}\int_{0}^{+\infty}e^{2i\frac{a}{b}\sigma ty^{2}}e^{-y^{2}}\ [T_{x}^{\nu,\mathbf{m}^{-1}}f]\left(\sqrt{4\sigma t}y\right) \ y^{2\nu+1}dy.\label{formula}
\end{eqnarray}
Since $\frac{2}{\Gamma(\nu+1)}\int_{0}^{+\infty}e^{-y^{2}}y^{2\nu+1}dy=1,$ the formula (\ref{formula}) implies
\begin{eqnarray}
S_{\nu}^{\mathbf{m}^{-1}}(t)f(x)-f(x)=a_{t}(x)+b_{t}(x)
\label{ref0}
\end{eqnarray}
where
\begin{eqnarray*}
a_{t}(x)&=&\frac{2}{\Gamma(\nu+1)}\int_{0}^{+\infty}e^{-y^{2}}\left(e^{2i\frac{a}{b}\sigma ty^{2}}-1\right)f(x)\ y^{2\nu+1}dy, \\
b_{t}(x)&=&\frac{2}{\Gamma(\nu+1)}\int_{0}^{+\infty}e^{2i\frac{a}{b}\sigma ty^{2}}e^{-y^{2}}\left([T_{x}^{\nu,\mathbf{m}^{-1}}f]\left(\sqrt{4\sigma t}y\right)-f(x)\right)\ y^{2\nu+1}dy.
\end{eqnarray*}
The relation (\ref{ref0}) implies that
\begin{eqnarray*}
\|S_{\nu}^{\mathbf{m}^{-1}}(t)f-f\|_{p,\nu}\leq \|a_{t}\|_{p,\nu}+\|b_{t}\|_{p,\nu}
\end{eqnarray*}
with
\begin{eqnarray*}
\|a_{t}\|_{p,\nu}\leq \frac{2}{\Gamma(\nu+1)}\left[\int_{0}^{+\infty}e^{-y^{2}}\left |e^{2i\frac{a}{b}\sigma ty^{2}}-1\right| y^{2\nu+1}dy\right]  \ \|f\|_{p,\nu}\longrightarrow 0 \ \mbox{as} \ t\rightarrow 0,
\end{eqnarray*}
and by Minkowski's inequality for integrals
\begin{eqnarray*}
\|b_{t}\|_{p,\nu}\leq \frac{2}{\Gamma(\nu+1)} \int_{0}^{+\infty}e^{-y^{2}}\left\|T_{\sqrt{4\sigma t}y}^{\nu,\mathbf{m}^{-1}}f-f\right\|_{p,\nu}y^{2\nu+1}dy\longrightarrow 0 \ \mbox{as} \ t\rightarrow 0.
\end{eqnarray*}
This is clear from the dominated convergence theorem since:
\begin{eqnarray*}
\left\|T_{\sqrt{4\sigma t}y}^{\nu,\mathbf{m}^{-1}}f-f\right\|_{p,\nu}&\leq& 2 \|f\|_{p,\nu}\quad (\mbox{by}\ (\ref{inequa})),\\
\ds\lim_{t\rightarrow 0}\|T_{\sqrt{4\sigma t}y}^{\nu,\mathbf{m}^{-1}}f-f\|_{p,\nu}&=&0\quad (\mbox{see Theorem \ref{theo13}}).
\end{eqnarray*}
Now suppose $f\in \mathcal{L}_{p,\nu}.$ If $\epsilon>0$ and since $\mathcal{C}_{\ast,c}(\mathbb{R})$ is dense in $\mathcal{L}_{p,\nu},$ there exists $g\in \mathcal{C}_{\ast,c}(\mathbb{R})$ such that $\|f-g\|_{p,\nu}\leq \frac{\epsilon}{3},$ so
\begin{eqnarray*}
\left\|S_{\nu}^{\mathbf{m}^{-1}}(t)f-f\right\|_{p,\nu}&\leq&\left\|S_{\nu}^{\mathbf{m}^{-1}}(t)(f-g)\right\|_{p,\nu}
+\left\|S_{\nu}^{\mathbf{m}^{-1}}(t)g-g\right\|_{p,\nu}+\left\|g-f\right\|_{p,\nu}\\
&\leq&2\left\|g-f\right\|_{p,\nu}+\left\|S_{\nu}^{\mathbf{m}^{-1}}(t)g-g\right\|_{p,\nu}\leq \frac{2\epsilon}{3}+\left\|S_{\nu}^{\mathbf{m}^{-1}}(t)g-g\right\|_{p,\nu}
\end{eqnarray*}
and $\left\|S_{\nu}^{\mathbf{m}^{-1}}(t)g-g\right\|_{p,\nu}\leq \frac{\epsilon}{3}$ if $t$ is sufficiently small.\\
$3)$  By Young's inequality (\ref{young})
\begin{eqnarray*}
\|u(t,.)\|_{r,\nu}&=&\|\mathcal{P}_{t}^{\mathbf{m}^{-1}}\underset{\nu,\mathbf{m}^{-1}}{*}f\|_{r,\nu}\\
&\leq& \|\mathcal{P}_{t}^{\mathbf{m}^{-1}} \ \|_{q,\nu}\ \|f\|_{p,\nu}.
\end{eqnarray*}
The desired result is therefore a consequence of
\begin{eqnarray*}
 \|\mathcal{P}_{t}^{\mathbf{m}^{-1}}\|_{q,\nu}&=&\frac{2}{\Gamma(\nu+1)(4\sigma t)^{(\nu+1)}}
\left(\ds\int_{0}^{+\infty}e^{-\frac{qy^{2}}{4\sigma t}}y^{2\nu+1}\ dy\right)^{1/q}\\
&=&\left(\frac{\Gamma(\nu+1)(4\sigma t)^{(\nu+1)}}{2}\right)^{\frac{1}{q}-1}\left(\frac{1}{q}\right)^{\frac{\nu+1}{q}}.
\end{eqnarray*}
$4)$Let $v(t,.)$ be another solution of (\ref{cauchy}) belonging in $\mathcal{L}_{r,\nu}$ with $r\in[1,+\infty)$. For $0\leq s<t<+\infty$
$$
\frac{d}{ds}\left[S_{\nu}^{\mathbf{m}^{-1}}(t-s)v(s,.)\right]
=S_{\nu}^{\mathbf{m}^{-1}}(t-s)\sigma\Delta_{\nu}^{\mathbf{m}^{-1}}v(s,.)-S_{\nu}^{\mathbf{m}^{-1}}(t-s)\sigma\Delta_{\nu}^{\mathbf{m}^{-1}}v(s,.)=0.
$$
Then $S_{\nu}^{\mathbf{m}^{-1}}(t-s)v(s,.)$ is independent of $s$; setting $s=0, s=t$ yeild
$$
v(t,.)=S_{\nu}^{\mathbf{m}^{-1}}(t)v(0,.)=S_{\nu}^{\mathbf{m}^{-1}}(t)f=u(t,.).
$$
$4)$ For $f\in \mathcal{S}_{\ast }(\mathbb{R})$ we have
\begin{eqnarray*}
u(t,x)=S_{\nu}^{\mathbf{m}^{-1}}(t)f(x)=\frac{2}{\Gamma(\nu+1)}\int_{0}^{+\infty}e^{-y^{2}}\ e^{2i\frac{a}{b}\sigma ty^{2}}[T_{x}^{\nu,\mathbf{m}^{-1}}f]\left(\sqrt{4\sigma t}y\right) \ y^{2\nu+1}dy.
\end{eqnarray*}
According to Theorem \ref{deri}, we can write:
\begin{eqnarray*}
e^{2i\frac{a}{b}\sigma ty^{2}}[T_{x}^{\nu,\mathbf{m}^{-1}}f][\sqrt{4\sigma t}y]&=&\ds\sum_{k=0}^{n}
\frac{\Gamma(\nu+1)}{\Gamma(k+\nu+1)k!}\ t^{k}\sigma^{k} \ y^{2k}\left(\Delta_{\nu}^{\mathbf{m}^{-1}}\right)^{k}f(x)\\
&+&\frac{t^{n+1}}{n!}\ds\int_{0}^{1}(1-s)^{n}\frac{d^{n+1}}{du^{n+1}}\left[e^{2i\frac{a}{b}\sigma uy^{2}}[T_{x}^{\nu,\mathbf{m}^{-1}}f][\sqrt{4\sigma u}y]\right](st)\ ds.
\end{eqnarray*}
Multiply both sides by $[2/\Gamma(\nu+1)] e^{-y^{2}}y^{2\nu+1}$ and integrate from $0$ to $+\infty,$ we get
\begin{eqnarray*}
u(t,x)&=&\ds\sum_{k=0}^{n}
 \frac{t^{k}}{k!}\sigma^{k}\left(\Delta_{\nu}^{\mathbf{m}^{-1}}\right)^{k}f(x)\left(\frac{2}{\Gamma(k+\nu+1)}\ds\int_{0}^{+\infty}e^{-y^{2}}y^{2(\nu+k)+1}\ dy\right)\\
&+&\frac{t^{n+1}}{n!}\left(\frac{2}{\Gamma(\nu+1)}
\ds\int_{0}^{+\infty}e^{-y^{2}}\left[\ds\int_{0}^{1}(1-s)^{n}\frac{d^{n+1}}{du^{n+1}}\left[e^{2i\frac{a}{b}\sigma uy^{2}}[T_{x}^{\nu,\mathbf{m}^{-1}}f][\sqrt{4\sigma u}y]\right](st)\ ds\right]y^{2\nu+1}dy\right).
\end{eqnarray*}
The desired result is therefore a consequence of $\frac{2}{\Gamma(k+\nu+1)}\int_{0}^{+\infty}e^{-y^{2}}y^{2(\nu+k)+1}\ dy=1.$
\end{Proof-}

\end{document}